\documentclass[12pt]{article}
\usepackage[backref=page,colorlinks,bookmarksopen=true]{hyperref}
\usepackage{amssymb}
\usepackage{latexsym}
\usepackage{amsmath}
\usepackage[all]{xy}
\newtheorem{theorem}{Theorem}[section]
\newtheorem{lemma}[theorem]{Lemma}
\newtheorem{proposition}[theorem]{Proposition}
\newtheorem{corollary}[theorem]{Corollary}
\newtheorem{definition}[theorem]{Definition}

\newtheorem{remark}[theorem]{Remark}

\newcommand{\Z}{\mathbb Z}

\newcommand{\C}{\mathcal{C}}

\newcommand{\M}{\mathcal{M}}

\begin{document}
\title{On $\Z/2\Z$-extensions of pointed fusion categories}
\author{Leonid Vainerman and Jean-Michel Vallin}
\date{1 October 2012}
\maketitle

\begin{abstract} We give a classification of $\mathbb Z/2\mathbb Z$-graded fusion
categories whose 0-component is a pointed fusion category.
A number of concrete examples is considered.  \end{abstract}

\tableofcontents

\footnote {AMS Subject Classification [2010]{: Primary 18D10, Secondary 16T05.}}

 \footnote{Keywords {: Fusion categories, extensions.}}
\newpage
\newenvironment{dm}{\hspace*{0,15in} {{\bf Proof.}}}{$\square$}
\begin{section}{Introduction}
\begin{subsection}{}A fusion category over an algebraically closed field $k$ of characteristic zero is a $k$-linear semisimple rigid tensor
category $\mathcal C$ with finitely many simple objects (the unit object $1$ is supposed to be simple) and finite dimensional
spaces of morphisms. Our main reference on fusion categories is \cite{ENO1}. Throughout this paper we work with $k=\mathbb C$
although many results hold for general $k$.
If there is no ambiguity, we use the same notation for an equivalence class and for its representative.

A fusion category $\mathcal C$ is said to be graded by a finite group $G$ if $\mathcal C=\oplus_{g\in G} c(g)$, where
$c(g)$ are full abelian subcategories of $\mathcal C$ such that $c(g)^*=c(g^{-1})$ and the tensor product maps $c(g)\times
c(h)$ to $c(gh)$, for all $g,h\in G$ (we call $\mathcal C$ a $G$-extension of its fusion subcategory $c(e)$, where $e$ is the
unit of $G$). A fusion category is said to be pointed if all its simple objects are invertible with respect to the
tensor product. Such a category is equivalent to the category $Vec_S^\omega$ whose simple objects are elements of a finite
group $S$ with tensor product $s\otimes t=st$, the unit object $1=1_S$, the duality $s^*={^*s}=s^{-1}$, and the associativity
isomorphisms defined by $\omega\in H^3(S,\mathbb C^\times)$.

Our aim is classification of $\mathbb Z/2\mathbb Z$-extensions $\mathcal C=c(0)\oplus c(1)$ of $c(0)=Vec_S^\omega$. Tambara and Yamagami
\cite{TY} showed that for any such an extension with $c(1)$ containing only one simple object, $S$ must be abelian
and equipped with a symmetric non-degenerate bicharacter $\chi$, and $\omega=1$. The general case is much more complicated. We
show that $S$ must contain a normal abelian subgroup $A$, so $S\cong A\underset{\rho}\rtimes S/A$ - a twisted semidirect product
with an action of $S/A$ on $A$ and $\rho\in Z^2(S/A,A)$ - see, for instance, \cite{B} (the simplest example: $S=\mathbb Z/4\mathbb Z=
\mathbb Z/2\mathbb Z\underset{\rho}\rtimes \mathbb Z/2\mathbb Z$ with trivial action and nontrivial $\rho$). Then, $A$ must be equipped with a symmetric
non-degenerate bicharacter $\chi$ and with an equivalence class of 2-cochains $\mu^r$ invariant with respect to $Aut(S)$ and
such that $\partial^2\mu^r=\omega$, modulo 2-cocycles on $A$ which can be extended to 2-cocycles on $S^{op}$, the opposite group to $S$.

The above extensions of 2-cocycles are classified by couples $(\psi,\nu)$, where map $\psi\in Z^1(S/A,Fun(A\times A,\mathbb C^\times))$ and $\nu\in C^2(S/A,\mathbb C^\times)$ satisfy some relations. If $S$ is a usual semidirect product, i.e., $\rho=1$, this was explained in \cite{Kar}; in this partial case we show that there is only one class $\mu^r$, so it
can be omitted. Finally, the Grothendieck rings of $\mathbb Z/2\mathbb Z$-extensions of $Vec_S^\omega$ are classified by couples
$(\varepsilon,\delta)$, where $\varepsilon\in Aut(S/A),\delta\in S/A$ are such that $\varepsilon^2=Ad(\delta)$ and $\varepsilon(\delta)=\delta$.

Thus, our main result - Theorem \ref{main} claims that $\mathbb Z/2\mathbb Z$-extensions of $Vec_S^\omega$ are classified,
up to equivalence, by collections $(A,\chi,\mu^r,\tau,\varepsilon,\delta,\psi,\nu)$, where  
$\tau=\pm|A|^{-1/2}$, like in \cite{TY}. It also describes their structure and allows to construct a number of new examples of fusion categories.

We have to mention that Liptrap \cite{Lip} earlier obtained some classification of $\mathbb Z/2\mathbb Z$-extensions of $Vec_S^\omega$
(analyzing, as in \cite{TY}, the solutions of the system of 16 pentagon equations for the associativity isomorphisms of
$\mathcal C$), but in terms that are difficult to apply to the
construction of concrete examples. On the other hand, we need new examples of fusion categories in order to construct new
families of finite index and finite depth $II_1$-subfactors as follows: 1) For a given a concrete fusion category
$\mathcal C$, the Hayashi's reconstruction theorem \cite{Ha} (see also \cite{SZ}) allows to construct a canonical weak Hopf
algebra (a quantum groupoid) $H$ \cite{BNSz}, \cite{NV} whose representation category is equivalent to $\mathcal C$.
2) Out of a given $H$, one can construct a subfactor whose bimodule category is equivalent to $\mathcal C$ and to compute its index,
principal and dual graphs, and the lattice of intermediate subfactors - see \cite{NV1}, \cite{NV2}. If $\mathcal C$
is Tambara-Yamagami category, this was done in \cite{M} and gave a family of subfactors of index $(n+\sqrt n)^2/d\
(n,d\in\mathbb N,\ d|n)$. We will describe in a separate paper a much larger family of subfactors coming from the fusion categories
constructed in the present paper.
\end{subsection}
\begin{subsection}{} The classification of $G$-extensions of fusion categories given in \cite{ENO2} implies, in the case when $G=\mathbb Z/2\mathbb Z$, that
any $\mathbb Z/2\mathbb Z$-extension
$\mathcal C=c(0)\oplus c(1)$ determines the following data (for all needed definitions see \cite{ENO2} and the references therein):

(1) A  group homomorphism $c:\mathbb Z/2\mathbb Z\to \pi_1=BrPic(c(0))\ (0\mapsto  \overline{c(0)}, 1\mapsto \overline{c(1)})$, where
the elements of the {\it Brauer-Picard group} $\pi_1$ are the equivalence classes of invertible $c(0)$-bimodule categories
and the operation in $\pi_1$ is the relative tensor product $\boxtimes_{c(0)}$. In fact, $c$ is defined by the choice of an
invertible $c(0)$-bimodule category $c(1)$ such that $c(1)\cong c(1)^{op}$.

(2) A collection of $c(0)$-bimodule equivalences $M_{g,h}:c(g)\boxtimes_{c(0)}c(h)\cong c(gh)$ $(g,h
\in\mathbb Z/2\mathbb Z)$ such that the following functors are isomorphic to $Id$:

$$
T_{f,g,h}:M_{fg,h}(M_{f,g}\boxtimes_{c(0)} Id_{c(h)})(Id_{c(f)}\boxtimes_{c(0)}M_{g,h}^{-1})
M_{f,gh}^{-1}:c(fgh)\to c(fgh).
$$

(3) Natural isomorphisms

$$
\alpha_{f,g,h}:M_{f,gh}(Id_{c(f)}\boxtimes_{c(0)}M_{g,h})\cong M_{fg,h}(M_{f,g}\boxtimes_{c(0)}Id_{c(h)})
$$

satisfying the pentagon equations

\begin{align*}
M_{f,gh,k}(id_{c(f)}&\boxtimes_{c(0)}\alpha_{g,h,k})\times \\
& \times \alpha_{f,gh,k}(Id_{c(f)}\boxtimes_{c(0)}M_{g,h}\boxtimes_{c(0)}Id_{c(k)})
M_{fgh,k}(\alpha_{f,g,h}\boxtimes_{c(0)} id_{c(k)})
\end{align*}
$$
=\alpha_{f,g,hk}(Id_{c(f)}\boxtimes_{c(0)}Id_{c(g)}\boxtimes_{c(0)}M_{h,k})\alpha_{fg,h,k}(M_{f,g}\boxtimes_{c(0)}
Id_{c(h)}\boxtimes_{c(0)}Id_{c(k)}).
$$

Vice versa, given $c(0)$, $c$ and  $M_{g,h}$ as above, the $c(0)$-bimodule category $\C=c(0)\oplus c(1)$ can be equipped with a
$c(0)$-bimodule tensor product which is associative if and only if certain cohomological obstruction $O_3(c)$ vanishes
or, equivalently, if and only if the functors $T_{f,g,h}$ are isomorphic to $Id$. If this is the case, the above tensor product
admits two possible choices of associativity isomorphisms satisfying the pentagon equations.
\end{subsection}
\begin{subsection}{} The paper is organized as follows:  Section 2 contains preliminary results on induction and extension of cocycles from
a subgroup of a finite group and also on invertible bimodule categories over $c(0)=Vec_S^\omega$. In Section 3 we
give a classification of $\mathbb Z/2\mathbb Z$-extensions of $c(0)=Vec_S^\omega$. In order to do this, we obtain the following intermediate results:

- Description of homomorphisms $c:\mathbb Z/2\mathbb Z\to BrPic(c(0))$,  or equivalently   invertible $c(0)$-bimodule categories $c(1)$ such that $c(1)\cong c(1)^{op}$.

- Classification of fusion rings of possible $\mathbb Z/2\mathbb Z$-extensions of $c(0)$.

- Explicit calculation of the $c(0)$-bimodule equivalences $M_{f,g}$.

- Explicit calculation of the functors $T_{f,g,h}$ which allows not only to deduce that they are isomorphic to $Id$ (i.e., that the cohomological obstruction $O_3(c)$ vanishes), but also to calculate explicitly the natural isomorphisms $\alpha_{f,g,h}$ satisfying the pentagon equations (there are exactly 2 choices of them).

Section 4 is devoted to examples: we compute the number of non-equivalent $\mathbb Z/2\mathbb Z$-extensions of $Vec_S^\omega$,
where $S$ is either an abelian group of order $2p$, or dihedral group $D_p\ (p$ is prime), or the alternate group $A_4$.

Note that our results imply the classification, up to (categorical) Morita equivalence, of $\mathbb Z/2\mathbb Z$-extensions of {\it group-theoretical categories}, i.e., Morita equivalent to pointed \cite{ENO1}. Indeed, \cite{ENO3}, Lemma 3.4 implies
that for any such an extension there is a Morita equivalent $\mathbb Z/2\mathbb Z$-extension of a pointed category.

{\bf Acknowledgements.} The authors are very grateful to Dmitri Nikshych for numerous discussions. In particular,
Proposition 2.7 is essentially due to him. They thank also Deepak Naidu for helpful comments.
\end{subsection}
\end{section}
\newpage
\begin{section}{Preliminaries}
\begin{subsection}{Some cohomological constructions}
\label{cohomology}

\begin{subsubsection}{\bf  Basic definitions,  induction of cocycles}
Let $C(S,P)=\{C^n(S,P)\}_{n\geq 0}$ be a cochain complex of a finite group $S$ with coefficients in a left or
right $S$-module $P$ \cite{B}. Namely, $C^n(S,P)=Fun(S^n,P)$ is a set of normalized $n$-cochains
(i.e., equal to $0$ if at least one of arguments equals to $1_S$), $S^n=S\times...\times S$ ($n$ factors).
If $P$ is a left $S$-module with an action $s\cdot p\ (s\in S,p\in P)$, then the coboundary operator $\partial^{n}:C^n(S,P)\to C^{n+1}(S,P)$ is defined by
$$
(\partial^{n} f)(s_1,.,s_{n},s_{n+1})=s_1\cdot f(s_2,..,s_{n+1})+\Sigma_{i=1}^n (-1)^i f(s_1,...,s_{i-1},s_i s_{i+1},..,s_{n+1})+
$$
\begin{equation} \label{lcobound}
+(-1)^{n+1}f(s_1,...,s_{n}).
\end{equation}
Similarly, if $P$ is a right $S$-module with an action $p\cdot s$, then the coboundary operator $\underline \partial^{n}:C^n(S,P)
\to C^{n+1}(S,P)$ is defined by
$$
(\underline\partial^{n} f)(s_1,.,s_{n},s_{n+1})=f(s_2,..,s_{n+1})+\Sigma_{i=1}^n (-1)^i f(s_1,...,s_{i-1},s_i s_{i+1},..,s_{n+1})+
$$
\begin{equation} \label{rcobound}
+(-1)^{n+1}f(s_1,...,s_{n})\cdot s_{n+1}.
\end{equation}
Let $Z^n(S,P)=Ker(\partial^{n})$ (resp., $\underline Z^n(S,P)=Ker(\underline\partial^{n})$) be the set
of $n$-cocycles, and also $B^n(S,P)=Im(\partial^{n-1})$ (resp., $\underline B^n(S,P)=Im(\underline\partial^{n-1})$)
the set of $n$-coboundaries, and $H^n(S,P)=$ $Z^n(S,P)/B^n(S,P)$ (resp., $\underline H^n(S,P)=\underline Z^n(S,P)/
\underline B^n(S,P)$) the $n$-th cohomology group of $S$ with coefficients in $P$.

Any left $G$-module is also a right $G$-module with the action $(s,M) \mapsto s^{-1}\cdot M$. Let $\sigma_n: C^{n}(S,P) \to C^{n}(S,P)$ be the map defined by $\sigma_n(f)(g_1,...,g_n) = -f(g_n^{-1},...g_1^{-1})$, we have easily $\sigma_{n+1}\partial^n\sigma_n = \underline  \partial^n$, so $\sigma_n$ is an isomorphism $\underline Z^n(S,P)\to Z^n(S,P) $ and passes to an isomorphism $\underline H^n(S,P)\to H^n(S,P)$.

Given a subgroup $A$ of $S$, we denote by $p:S\to S/A$ the usual surjection $p(s)=sA$, for all $s\in S$, and ${\bf 1}=p(1_S)$. Let us
choose a representative $u(M)$ in any coset $M\in S/A$, in particular, $u({\bf 1})=1_S$. $S$ acts on $S/A$ via $s\cdot M=p(su(M))$
and also on the set $\{u(M)|M\in S/A\}$ via $s\cdot u(M)=u(s\cdot M)$.
Then, for all $s\in S,M\in S/A$, there exists an element $\kappa_{M,s}\in A$ such that $su(M)=u(s\cdot M)\kappa_{M,s}$.
One can check that $\kappa_{M,s_1s_2}=\kappa_{s_2\cdot M,s_1}\kappa_{M,s_2}$.

Let $C=Fun(S/A,\mathbb C^\times)$ be the coinduced right $S$-module with the natural action $f(M)\cdot s=f(s\cdot M)$,
for all $s\in S,M\in S/A$ (here $\mathbb C^\times$ is viewed as a trivial right $S$-module). By Shapiro's lemma (see \cite{B})
the groups $\underline H^n(S,C)$ and $\underline H^n(A,\mathbb C^\times)$ are isomorphic. Explicitly, \cite{Naidu}, Lemmas 2.1
and 2.2 show that for $n=1$ this isomorphism is induced by the maps
$$
\varphi_1:\underline Z^1(A,\mathbb C^\times)\to\underline  Z^1(S,C): (\varphi_1(\rho)(s))(M)=\rho(\kappa_{M,s}),
$$
\begin{equation} \label{shapiro1}
\varphi^{-1}_1:\underline Z^1(S,C)\to\underline  Z^1(A,\mathbb C^\times):\ \varphi^{-1}_1(\beta)(a)=\beta(a)({\bf 1}),
\end{equation}
and for $n=2$, respectively, by the maps
$$
\varphi_2:\underline Z^2(A,\mathbb C^\times)\to\underline Z^2(S,C): (\varphi_2(\mu)(s_1,s_2))(M)=\mu(\kappa_{s_2\cdot M,s_1},\kappa_{M,s_2}),
$$
\begin{equation} \label{shapiro2}
\varphi^{-1}_2:\underline Z^2(S,C)\to\underline Z^2(A,\mathbb C^\times):\ \varphi^{-1}_2(\gamma)(a_1,a_2)=\gamma(a_1,a_2)({\bf 1}).
\end{equation}
\end{subsubsection}
\begin{subsubsection}{\bf   Extension of 2-cocycles}
If $A\triangleleft S$ is abelian, then $S$ is isomorphic (see \cite{B}, IV, 3)) to the twisted semi-direct product
$A\underset \rho\rtimes T$ defined by:
$$
(a,t)\cdot(a',t')= ((a\cdot\ ^t a')\rho(t,t'),tt')\ \ \forall a,a'\in A,\ t,t'\in T,
$$
where $T = S/A$ acts on $A$ by inner automorphisms (i.e., $^t a = u(t)au(t)^{-1},$ $u(t)\in S$), $\rho\in Z^2(T,A)$ is given by $\rho(t,t') = u(t)u(t')u(tt')^{-1}$, the isomorphism between $A\underset \rho \rtimes T$ and $S$ is defined by $(a,t) \mapsto au(t)$. Moreover, the map
$(a,t)^{op} \mapsto (^{t^{-1}} a,t)$ is an isomorphism between $S^{op}$ and $A\underset {\rho^{op}}\rtimes T^{op}$, where $\rho^{op}(t,t') = ^{(t't)^{-1}}\rho (t',t)=u(t't)^{-1}u(t')u(t)\in Z^2(T^{op},A)$ and the action of $T^{op}$ on $A$ is given by $^{t^{op}} a =^{t^{-1}} a$.

In the case of usual semidirect product (i.e., $\rho=1$) Karpilovski \cite{Kar} explained how to extend $\sigma\in H^2(A,
\mathbb C^\times)$ to $\mu\in H^2(S,\mathbb C^\times)$. We generalize this construction to twisted semidirect products using essentially the same arguments (however, slightly more complicated because of the presence of $\rho$). Like in \cite{Kar},
Lemma 2.2.3, one can show that:

1) Any $\mu\in Z^2(S,\mathbb C^\times)$ is cohomologous to $\mu'$ that is {\it normal}, i.e., $\mu'((a,e),$ $(e,t))=1$, for all $a\in A,t\in T$, and has the same restriction $\mu_{T,T}$ on $(e,T)\times (e,T)$. Note that $(e,T)$ is not a subgroup of $S$, in general, because the products in $S$ and in $T$ are related by the formula $(e,t)\cdot_S(e,t)'=(\rho(t,t'),tt')$.

2) Any normal $\mu\in Z^2(S,\mathbb C^\times)$ is completely determined by its restrictions $\mu_{T,T}$, $\mu_{A,A}:=\mu|_{(A,e)\times(A,e)}$, and $\mu_{T,A}:=\mu|_{(e,T)\times(A,e)}$ by the following formula in which we identify $(a,t)$ with $at$ and $(a',t')$ with $a't' (\forall
a,a'\in A,t,t'\in T)$:
\begin{equation} \label{extcoc}
\mu(at,a't') = \mu_{T,T}(t,t')\mu_{T,A}(t,a')\mu_{A,A}(a, ^t a')\mu_{AA}(a(^ta'),\rho(t,t')).
\end{equation}
\begin{proposition} \label{lift}
A 2-cocycle $\sigma\in Z^2(A,\mathbb C^\times)$ can be extended to a normal 2-cocycle $\mu\in Z^2(S,\mathbb C^\times)$
if and only if:

1) $\sigma$ is cohomologically $S/A$-invariant, i.e., there exists $\mu_{T,A}\in C^1(S/A,$ $Fun(A,\mathbb C^\times))$ such that
$\displaystyle \frac{\sigma}{^t\sigma}=\partial^1_A\mu_{T,A}(t,\cdot)$, where $t\in T,^t\sigma(a,b):=\sigma(^t a,^t b)$.

2) $\mu_{T,A}\in Z^1(T,Fun(A,\mathbb C^\times))$, where $^{t'}\mu_{T,A}(t,a):=\mu_{T,A}(t,^{t'}a),\ \forall\ t,t'\in T,a\in A$.

3) The 3-cochain on $T$
$$
\zeta(t,t',t"):=\mu_{T,A}(t,\rho(t',t''))\frac{\sigma(^t\rho(t',t''),\rho(t,t't''))}{\sigma(\rho(t,t'),\rho(tt',t''))}
$$
is a 3-coboundary, i.e., can be presented as $\partial^2\mu_{T,T}$ for some $\mu_{T,T}\in C^2(T,\mathbb C^\times)$.
\end{proposition}
\begin{dm}
1) Condition 1) is necessary due to \cite{Kar}, Proposition 1.5.8 (iii).

2) If a normal extension $\mu$ of $\sigma$ exists, write for it the 2-cocycle equality restricted to $(e,T)\times (e,T)\times(e,T)$, in terms of the product in $S$. Passing then to the product in $T$ and using (\ref{extcoc}),
we get the necessity of condition 3).

3) By direct calculations, exactly like in \cite{Kar}, Lemma 2.2.4, one has:
$$
\mu_{T,A}(tt',a)=\mu_{T,A}(t',a)\mu_{T,A}(t,^{t'}a)
$$
with the product $tt'$ in $T$. This means that $\mu_{T,A}\in Z^1(T,Fun(A,\mathbb C^\times))$.

Vice versa, it is straightforward to check that the relation (\ref{extcoc}) defines, under conditions 1),2),3), a normal extension of $\sigma$.
\hfill\end{dm}
\begin{remark} \label{extcohom}

1) In a slightly more general situation when 2-cocycles $\mu$ and $\sigma$ are replaced by 2-cochains such that
$\partial^2\mu=\partial^2 h$ for some $h\in C^2(S,\mathbb C^\times)$ and $\partial^2\sigma=\partial^2 h|_{A\times A\times A}$,
respectively, one can use Proposition \ref{lift} applied to 2-cocycle $\displaystyle \frac{\sigma}{h_{\mid A\times A}}$ which describes when this 2-cocycle can be extended
to a normal 2-cocycle $\displaystyle \frac{\mu}{h}$ or, equivalently, when $\sigma$ can be extended to $\mu$.

2) Let $\sigma=\partial^1 \eta$ be a cohomologically $S/A$-invariant 2-coboundary on $A$ and $\mu$ its normal extension described by $\mu_{T,A}$
and $\mu_{T,T}$ as above. Then one can deduce from (\ref{extcoc}) that $\mu$ is a 2-coboundary if and only if
$\mu_{T,A}=\partial^0 \tilde\eta$ and $\mu_{T,T}= \displaystyle \frac{\partial^1 f_T}{\eta\circ\rho}$, where $f_T\in C^1(T,\mathbb C^\times)$ and $\tilde\eta
(t,\cdot)\equiv \eta$ is a constant function.
\end{remark}
\end{subsubsection}
\end{subsection}
\begin{subsection}{Invertible bimodule categories over $c(0)=Vec_{S}^\omega$}
\label{inv}
\begin{subsubsection}{\bf  Indecomposable left $c(0)$-module categories} These categories are indexed (see \cite{O2}) by conjugacy classes of pairs $(A,\mu)$,
where $A$ is a subgroup of $S$ and a 2-cochain $\mu\in C^2(A,\mathbb C^\times)$ satisfies $\partial^2\mu=\omega\vert_{A\times A\times A}$
(so $\mu|_{A\times A\times A}=1$ in $H^3(A,\mathbb C^\times)$). A 2-cochain $\mu$ satisfying the relation $\partial^2\mu=\omega$ will be
called {\it $\omega$-2-cocycle.}

If $\mathcal M(A,\mu)$ is such a category, then group $S$ acts transitively on the left on the set $Irr(\mathcal M(A,\mu))=S/A$ of its simple
objects; its associativity isomorphisms are defined by a 2-cochain $\tilde\mu(s,t,u)\in C^2(S,C)$ induced from $\mu$ and such that
\begin{equation} \label{lmod}
\tilde\mu(t,u,M)\tilde\mu(s,tu,M)=\omega(s,t,u)\tilde\mu(s,t,u\cdot M)\tilde\mu(st,u,M).
\end{equation}
In its turn, $\mu(a,b)=\tilde\mu(a,b,\bf{1})$ for all $a,b\in A$.
\begin{remark} \label{ind}
If $\omega=1$, the induction above is given explicitly by the map $\varphi_2$ from subsection \ref{cohomology}. In general, there is no
a canonical way for such an induction. One can proceed as follows. Fix $\mu_0\in C^2(A,\mathbb C^\times)$ such that $\partial^2\mu_0=
\omega|_{A\times A\times A}$ and $\tilde\mu_0\in C^2(S,C)$ such that $\partial^2\tilde\mu_0=\omega$. Then, for any other $\mu\in C^2(A,
\mathbb C^\times)$ such that $\partial^2\mu=\omega|_{A\times A\times A}$, we can put $\tilde\mu=\tilde\mu_0 \varphi_2(\mu/\mu_0)$.
\end{remark}
Two pairs, $(A,\mu)$ and $(A',\mu')$, give rise to equivalent $c(0)$-module categories if and only if $A'=sAs^{-1}$ for
some $s\in S$ and $\mu$ is cohomologous to the $s$-conjugate $(\mu')^s$ of $\mu'$, that is, they differ by a 2-coboundary. Let
$$
\Omega_{A,\omega}:=\text{equivalence\ classes\ of}\ \{\mu\in C^2(A,\mathbb C^\times)\vert{\partial^2}\mu=\omega|_{A\times A\times A}\}.
$$
Note that $\Omega_{A,\omega}$ is a torsor over $H^2(A,\mathbb C^\times)$, in particular, $\Omega_{A,1}=H^2(A,\mathbb C^\times)$. For any $s\in S$ and $\mu\in\Omega_{A,\omega}$, let us define
$$
\mu\triangleleft s:=\mu^s\times \Upsilon_s\vert_{A\times A},\hskip 1cm\text{where}
$$
$$
\Upsilon_s(t,u):=\frac{\omega(sts^{-1},sus^{-1},s)\omega(s,t,u)}{\omega(sts^{-1},s,u)},\ \mu^s(t,u):=\mu(sts^{-1},sus^{-1}),
$$
for all $\ s,t,u\in S$. Let $(\Omega_{A,\omega})^S$ be a set of $S$-invariant elements of $\Omega_{A,\omega}$, i.e.,
$$
(\Omega_{A,\omega})^S:=\{\mu\in\Omega_{A,\omega}\vert\mu^s\times \Upsilon_s\vert_{A\times A}=\mu\ \text{in}\ H^2(A,\mathbb C^\times),\ \text
{for\ all}\ s\in S\}.
$$
If $\omega=1$, we have a usual definition of a class of $S$-invariant 2-cocycles on $A$.
\end{subsubsection}
\begin{subsubsection}
{\bf Invertible $c(0)$-bimodule categories} By definition, a $c(0)$-bimodule category $\mathcal M$ is a left $Vec_{S\times S^{op}}^{\omega\otimes\omega^{op}}$-module
category, where $S^{op}$ is the group opposite to $S$ and also $\omega^{op}(s^{op},t^{op},u^{op})=\omega^{-1}(s^{-1},t^{-1},u^{-1})$, for all
$s,t,u\in S$. The action of $(s,t^{op})\in S\times S^{op}$ on $M\in Irr(\mathcal M)$ defines left and right actions of $S$ by
$(s,t^{op})\cdot M:=(s\cdot M)\cdot t$,
so $\mathcal M$ can be viewed as both left and right $c(0)$-module category. Note that right indecomposable $c(0)$-module categories
are also parameterized by the classes of equivalence of pairs $(A,\mu)$, their associativity isomorphisms
are defined by 2-cochains $\tilde\mu(s,t,M)\in C^2(S,C)$ induced from $\mu$ and satisfying
\begin{equation} \label{rmod}
\tilde\mu(M\cdot s,t,u))\tilde\mu(M,s,tu)=\omega(s,t,u)\tilde\mu(M,s,t)\tilde\mu(M,st,u).
\end{equation}

If $\mathcal M$ is {\it invertible} (i.e., $\mathcal M^{op}\boxtimes_{c(0)}\mathcal M\cong\mathcal M\boxtimes_{c(0)}\mathcal M^{op}
\cong c(0)$, where $\mathcal M^{op}$ is the $c(0)$-bimodule category opposite to $\mathcal M$ and $\boxtimes_{c(0)}$
is the relative tensor product - see \cite{ENO2}), then it is indecomposable as both left and right $c(0)$-module category (see \cite{ENO2},
Corollary 4.4), so it is indecomposable as a left $Vec_{S\times S^{op}}^{\omega\otimes\omega^{op}}$-module category. Thus, it is of
the form $\mathcal M(L,\mu)$, where $L<S\times S^{op}$ and $\mu\in C^2(L,\mathbb C^\times)$ satisfies $\partial^2\mu=
(\omega\otimes\omega^{op})\vert_{L\times L\times L}$. Its associativity isomorphisms are defined by a 2-cochain $\tilde\mu\in C^2(S\times S^{op},C)$
induced from $\mu$.

Since $S$ acts transitively on $(S\times S^{op})/L$ on both sides, we have
$$
(S\times \{e\})L=(\{e\}\times S^{op})L=S\times S^{op}.
$$
Let $A_1$ be the subgroup of $S$ such that $L\cap(S\times\{e\})=A_1\times\{e\}$ and $A_2$ be the  subgroup of $S^{op}$ such that $L\cap(\{e\}\times S^{op})=\{e\}\times A_2$, and let us denote $\mu^l=\mu|_{(A_1,e)\times(A_1,e)},\mu^r=\mu|_{(e,A_2)\times(e, A_2)}$. Then $\mathcal M(L,\mu)$ viewed
as a left (resp., right) $c(0)$-module category is equivalent to $\mathcal M(A_1,\mu^l)$ (resp., to $\mathcal M(A_2,\mu^r)$).
\begin{lemma}
\label{rightaction}
The maps $f_1: (S\times S^{op})/L\to S/A_1$ and $f_2: (S\times S^{op})/L\to S^{op}/A_2$ defined by $M\mapsto p_1(M\cap(S\times \{e\}))$ and $M\mapsto p_2(M\cap (\{e\}\times S^{op}))$, respectively, where $p_1:(s,t^{op})\mapsto s, \ \ p_2:(s,t^{op})\mapsto t^{op}$, for all $s,t\in S$, are well defined bijections between $Irr(\mathcal M)$ and $S/A_1$ (resp., $  S^{op}/A_2$). Let us define $Inv: S^{op}/A_2 \to A_2 \setminus S$ by $s^{op}A_2 \mapsto A_2s^{-1}$, then the composition $f:=Inv\circ f_2\circ f_1^{-1}$ is a natural bijection:
$$
f : S/A_1 \to A_2 \setminus S
$$
such that:
\begin{equation}
\label{parameterized L}
L =\{(s,t^{op}) \in S\times S^{op} / f(sA_1)  = A_2t  \}.
\end{equation}
We will denote $L$ in (\ref{parameterized L}) by $L(A_1,\, A_2,\, f)$ and, if $A_1=A_2=A$, - by $L(A,f)$.

The  left action of $S$ (resp., $S^{op}$) on $Irr(\mathcal M)$ identified with  $S/A_1$ (resp., $S^{op}/A_2$), is given by the multiplication of left classes. The left action of $S^{op}$ on $Irr(\mathcal M)$ identified with $S/A_1$, is the right action of $S$ on $S/A_1$ given by:
$$
(xA_1)\cdot s = xf^{-1}(A_2s^{-1}).
$$
\end{lemma}
\begin{dm}
As $(S \times \{e\})L = S \times S^{op}$, then for any $s,t \in S$ there  exists $\sigma \in S$ such that $(\sigma,e)\in (s,t^{op})L$, so $p_1(M\cap(S\times\{e\}))$ is a not empty class in $S/A_1$, for any $M \in Irr(\mathcal M)$. Thus, the map $f_1$ above
is a bijection of $Irr(\mathcal M)$ and $S/A_1$, and similarly for $f_2$. One can see that the natural left action of $S$ (resp., $S^{op}$) on $(S\times S^{op})/L$ becomes the natural left action of $S$ (resp., $S^{op}$) on $S/A_1$ (resp., $S^{op}/A_2$). Moreover, for all $s,t \in S, f(sA_1) =  A_2t$ means $(s,e)L = (e,(t^{op})^{-1})L$ which is equivalent to $(s,t^{op})\in L$. So,  $L=\{(s,t^{op})\in S\times S^{op}\vert f(sA_1) = A_2t \}$. For any $M\in Irr(\mathcal M)$, there exists $x \in S$ such that $M = (x,e)L$ and, after identification with $S/A_1$, $M = xA_1$. One also has, for any $s\in S$:
$(e,s^{op})M = (x,s^{op})L$; if $y \in S$ is such that $(x,s^{op})L = (y,e)L$, this means $(x^{-1}y,(s^{op})^{-1})\in L$, hence $yA_1= xf^{-1}(A_2s^{-1})$.
\hfill\end{dm}

Let us describe $\mathcal M(L,\mu)$ in other terms - see \cite{Green1}, 2.1.1. Using notations $\tilde\mu^l(x,y,M):= \tilde \mu((x,e),(y,e),M),\ \tilde \mu^r(M,x,y) :=\tilde \mu((e,y^{op}),(e,x^{op}), M)$, and finally also let us define $\tilde\chi(x,M,y):=\tilde \mu((x,e),(e,y^{op}), M)$ and the decomposition
\begin{equation}\label{relbimod}
\tilde \mu((x_1,x^{op}_2),(y_1,y^{op}_2),M)=\tilde\chi(x_1,y_1\cdot M,y_2)\tilde\mu^l(x_1,y_1,M)\tilde\mu^r((x_1y_1)\cdot M,y_2,x_2),
\end{equation}
(see \cite{Green2}, p. 27), one can check that
$\underline{\partial}^2\tilde\mu^l=\omega, \partial^2\tilde\mu^r=\omega$ and that the following compatibility conditions
hold:
\begin{equation}\label{leftbimod}
\tilde\mu^l(x,y,M\cdot z)\tilde\chi(xy,M,z)=\tilde\chi(y,M,z)\tilde\chi(x,y\cdot M,z)\tilde\mu^l(x,y,M),
\end{equation}
\begin{equation}\label{rightbimod}
\tilde\mu^r(M,y,z)\tilde\chi(x,M,yz)=\tilde\chi(x,M\cdot y,z)\tilde\chi(x,M,y)\tilde\mu^r(x\cdot M,y,z).
\end{equation}
Here $\tilde\mu^l$ and $\tilde\mu^r$ define left and right $Vec_S^\omega$-module category structures
on $\mathcal M(L,\mu)$, respectively, $x,x_1,x_2,y,y_1,y_2,z\in S,,M\in Irr(\mathcal M(L,\mu))$.
\begin{lemma}
\label{toutecrire} If $\mathcal M(L,\mu)$ is an invertible $c(0)$-bimodule category, then:
(i) $A_1$ and $A_2$ are normal abelian subgroups equipped with $S$-invariant $\omega$-2-cocycles
 $\mu^l$ and $\mu^r$, respectively.
(ii) The map $f: S/A_1 \xrightarrow{\sim} S/A_2$ in Lemma \ref{rightaction} is a group anti-isomorphism.
\end{lemma}

\begin{dm}
\cite{ENO2}, Definition 4.1 and Proposition 3.5 imply that the dual category of $c(0)$ with respect to $\mathcal M(L,\mu)$
viewed as a right $c(0)$-module category (and so equivalent to $\mathcal M(A_2,\mu^r)$), is equivalent to $c(0)$ itself which
is pointed. But due to \cite{Naidu}, Theorem 3.4, this is possible if and only if the pair $(A_2,\mu^r)$ satisfies conditions (i).
Similarly for $(A_1,\mu^l)$.

(ii) Follows from (i) and from Lemma \ref{rightaction}.\hfill\end{dm}

In what follows we denote by $\chi(\cdot,\cdot)$ the bicharacter $\tilde\chi(\cdot,{\bf 1},\cdot)|_{A_1\times A_2}$.

{\bf 3.} The additive endofunctor $L(s)$ of $\mathcal M(L,\mu)$ defined by left multiplication by $(s,e)$, where $s\in S$, is isomorphic to the
identity if and only if $s\in A_1$. Due to (\ref{rightbimod}), its right $c(0)$-module functor structure can be defined by
$$
\tilde\chi^{-1}(s,M,x)id_{s\cdot M\cdot x}:
L(s)(M\cdot x) \xrightarrow{\sim} (L(s)(M))\cdot x,\qquad s,x\in S.
$$
Similarly, the additive endofunctor $R(s)$ of $\mathcal M(L,\mu)$ defined by left multiplication by $(e,s^{op})$,
has a structure of a left $c(0)$-module functor:
$$
\tilde\chi(x,M,s)id_{x\cdot M\cdot s}:
R(s)(x\cdot M) \xrightarrow{\sim} x\cdot(R(s)(M)),\qquad x,s\in S.
$$

Let us summarize these observations:

\begin{lemma}
\label{when id}
$L(s)$ is equivalent to the identity as a right $c(0)$-module autoequivalence of
$\mathcal M(L,\mu)$ if and only if $s\in A_1$ and $a_{s}(\cdot) =1$ on $A_2$,
where the group homomorphism $a: A_1 \to \widehat{A_2}$ is defined by
\begin{equation}
\label{a alt}
a: s \mapsto a_s(x):=\chi^{-1}(s,x),\, x\in A_2.
\end{equation}
\end{lemma}

Similar statement is valid for $R(s)$.

\begin{proposition}
\label{M(L,mu) inv}
$\M(L,\mu)$ is invertible if and only if:
\begin{enumerate}
\item[(i)] The conditions of Lemma \ref{toutecrire} hold.
\item[(ii)] There is $\mu'$ cohomologous to $\mu$ such that the well  defined $S$-invariant bicharacter
$\chi:=\mu'|_{(A_1,e)\times(e,A_2)}: A_1 \times A_2 \to \mathbb C^\times$ is non-degenerate.
\end{enumerate}
\end{proposition}

\begin{dm}
The only thing to prove is (ii). \cite{ENO2}, Prop. 4.2 claims that $\M=\M(L,\mu)$ is invertible if and only if the functors
$c(0) \to Fun_{c(0)}(\M,\, \M): s \mapsto L(s)$ (respectively, $c(0) \to Fun(\M,\, \M)_{c(0)}: s \mapsto R(s)$) are equivalences.
Since those functors are tensor, the latter condition is equivalent to $L(s) \not\cong id_\M$ as a right $c(0)$-module functor
(respectively, $R(s)\not\cong id_\M$ as a left $c(0)$-module functor), for all $s\neq e$.
Due to Lemma~\ref{when id}, these conditions hold if and only if the group homomorphisms \eqref{a alt} and
\begin{equation}
\label{rho2}
A_2 \to \widehat{A_1} :  x \mapsto a'_x,
\mbox{where } a'_x(s) :=\chi(s,x),\, s\in A_1,
\end{equation}
are injective. This is equivalent to $\chi$ being non-degenerate on $A_1\times A_2$.
\hfill\end{dm}
\end{subsubsection}
\end{subsection}
\end{section}
\begin{section}{$\mathbb Z/2\mathbb Z$-extensions of $c(0)=Vec_S^\omega$}

\begin{subsection} {Group homomorphisms $c:\mathbb Z/2\mathbb Z\to BrPic(Vec_S^\omega)$ and fusion rings}
\begin{subsubsection}{\bf Homomorphisms $c:\mathbb Z/2\mathbb Z\to BrPic(Vec_S^\omega)$ }
 Such a   group homomorphism $c$   is defined by an invertible $c(0)$-bimodule category $\mathcal M(L,\mu)$ which
is equivalent to its opposite. The last one, $\mathcal M^{op}(L,\mu)$, has the same simple objects as $\mathcal M(L,\mu)$ on which $S\times S^{op}$ acts as $(s,t^{op})\cdot_{op} M=(s,t^{op})^\vee M$, where $\vee:(s,t^{op})\mapsto (t^{-1},(s^{-1})^{op})$ is canonical involutive automorphism of $S\times S^{op}$ (see \cite{ENO2}). The corresponding subgroup of $S\times S^{op}$ is $L^{op}=\vee(L)$ and $\mu_{op}=\mu(\vee\times\vee)^{-1}$. Passing to the quotient we have a canonical bijection $\overline \vee:(S\times S^{op})/L\mapsto(S\times S^{op})/L^{op}$ which transforms the opposite action of $S\times S^{op}$ on $(S\times S^{op})/L$ into the left natural action $\cdot_\vee$ of $S\times S^{op}$ on $(S\times S^{op})/L^{op}$, and one can check that $\mathcal M^{op}(L,\mu)\cong\mathcal M^{op}(L^{op},\mu_{op})$.
In particular, $\mu_{op}^l=(\mu^r\circ\vee)^{-1},\mu_{op}^r=(\mu^l\circ\vee)^{-1},\chi_{op}=\chi\circ\vee$.

In terms of Lemma \ref{toutecrire}, $L(A_1,\, A_2,\, f)^{op}= L(A_2,\, A_1,\,$ $f^{-1})$, so $\M(L,\mu)$ $\cong \M(L,\mu)^{op}$
if and only if subgroups $L(A_1,\, A_2,\, f)$ and $L(A_2,\, A_1,\,$ $f^{-1})$ are conjugate and $\mu_{op} $ is cohomologous
to the conjugate of $\mu$. In particular, this implies $A_1=A_2=A$, so $f$ is an anti-isomorphism of $S/A$.
\begin{lemma} \label{a86}
If $f$ and $f'$ are anti-isomorphisms of $S/A$, then $L=L(A,f)$ and $L'=L(A,f')$ are conjugate subgroups of $S\times S^{op}$
if and only if $f=Ad(xA)\circ f'$ for some $x\in S$. For such an $x$, we have $L'=Ad(1,x)L$ in $S\times S^{op}$.
In particular,  $L$ and $L^{op}$ are conjugate if and only if
\begin{equation}
\label{prepa}
f^2 = Ad(xA)\hskip 1cm\text{for\ some}\ \ \ x \in S.
\end{equation}
For such an $x$, we have $L^{op} = Ad(x,1)L= Ad(1,x)L$ in $S \times S^{op}$.
\end{lemma}

\begin{dm}
The sufficiency is clear. Let $(h,k)\in S\times S^{op}$ be such that $L'= Ad(h,k)L$. If $(s',t')\in L'$, then $f'(s'A)=t'A$,
but $Ad(h,k)^{-1}(s',t')\in L$, so $f((hA)(s'A)(hA)^{-1}) =(kA)(t'A)(k^{-1}A)$, from where $f=Ad(xA)\circ f'$ for
$xA=f(hA)^{-1}(kA)$.

For any $x \in S$ such that $f= Ad(xA)\circ f'$ and any $(s',t')\in L'$, one has in $S\times S^{op}$: $Ad(1,x^{-1})
(s',t') =(s',xt'x^{-1})$, but $(xA)(t'A)(x^{-1}A) =Ad(xA)\circ f'(s'A) = f(s'A)$, so $Ad(1,x^{-1})(s',t')\in L$. The proof of the
remaining statements is similar.\hfill\end{dm}

\begin{remark} \label{cross2}
If $\mathcal M(L_1,\mu_1)$ and $\mathcal M(L_0,\mu_0)$ are equivalent invertible $c(0)$-bimodule categories, then $L_1(A_1,f_1)$ and
$L_0(A_0,f_0)$ are conjugate and $\mu_1$ is cohomologous to a conjugate of $\mu_0$, so $A_1=A_0=A$ and $f_1=Ad(xA)\circ f_0$ for some
$x\in S$. In terms of the triples $(\tilde\mu^l_1,\tilde\mu^r_1,\tilde\chi_1)$ and $(\tilde\mu^l_0,\tilde\mu^r_0,\tilde\chi_0)$, this
equivalence is a bijection $F$ of $S/A$ equipped with left and right $\mathcal D$-module functor structures $f^l(s,M):F(sM)\mapsto
sF(M), f^r(M,s):F(Ms)\mapsto F(M)s$ such that, for all
$s,t\in S,M\in Irr(\mathcal M(L_1,\mu_1))$:
$$
\tilde\mu^l_1(s,t,F(M))f^l(st,M)=f^l(t,M)f^l(s,t\cdot M)\tilde\mu^l_0(s,t,M),
$$
$$
\tilde\mu^r_1(F(M),s,t)f^r(M,st)=f^r(M,s)f^r(M\cdot s,t)\tilde\mu^r_0(M,s,t),
$$
$$
\tilde\chi_1(s,F(M),t)f^l(s,M)f^r(s\cdot M,t)=f^r(M,t)f^l(s,M\cdot t)\tilde\chi_0(s,M,t)
$$
(see \cite{Green1}, Remark 2.14). This implies that $\mu^l_1$ and $\mu^r_1$ are cohomologous, respectively, to $\mu^l_0$ and $\mu^r_0$,
and that $\chi_1=\chi_0$.
\end{remark}

As $S\cong A\underset{\rho}\rtimes T,\ T=S/A$ (see \cite{B}, IV, 3), we have the following
\begin{lemma}
\label{complot}
Identifying $f$ with an anti-automorphism $\phi$ of $T$ such that $\phi^2$ is inner, we have $L(A,f)\cong (A\times A)
\underset{\tilde\rho}\rtimes T$ with the action of $T$ on $A\times A$ given by $^{\overline t} (a,a')=
(^t a,^{\phi(t^{-1})} a')$ and the 2-cocycle $\overline\rho(t,t') = (\rho(t,t'),^{(t't)^{-1}}\rho (t',t))$.

The map $\tilde \phi$: $((x,t),(y,s))\mapsto s^{-1}\phi(t)$ factors through $L$ to a bijection $(S\times S^{op})/L\to T$. If
$T$ is abelian, $L(A,f)\triangleleft (S\times S^{op})$ and $Irr(\mathcal M(L,\mu))$ is a group isomorphic to $T$ via the canonical
decomposition of $\tilde \phi$.
\end{lemma}

\begin{dm} Obviously,
$$
L(A,f) = \{((a,t),(b,\phi(t))/a,b \in A, t \in T\},
$$
and the map $((a,t),(b,\phi(t))\mapsto (a,\phi(t)^{-1}b,t)$ is an isomorphism $L(A,f)\to (A\times A)\underset {\overline\rho}\rtimes T$
with the above action of $T$ on $A\times A$ and 2-cocycle $\overline\rho$.
Let $ x= ((a,t),(b,s))$ and $y=((a',t'),(b',s'))$, then $x^{-1}y\in L(A,f)$, means that $\phi(t^{-1}t') = s's^{-1}$ or
$s'^{-1}\phi(t') =  s^{-1}\phi(t)$, the remaining statements follow.\hfill\end{dm}

\begin{corollary} \label{triple}
1) Given a triple $(\mu^l,\mu^r,\chi)$ as above, let us construct the $\omega\otimes\omega^{op}$-2-cocycle on $A\times A$:
$$
\mu_0((a,b),(a',b')):=\mu^l(a,a')\mu^r(b,b')\chi(a,b').
$$
In fact, $\mu_0=\mu|_{A\times A}$. Let $h\in C^2(L,\mathbb C^\times)$ be such that $(\omega\otimes\omega^{op})|_{L\times L\times L}=\partial^2 h$. Then
there exists an $\omega\otimes\omega^{op}$-2-cocycle $\mu$ on $L$ such that $\mu^l=
\mu|_{(A,e)\times(A,e)},\mu^r=\mu|_{(e,A)\times(e,A)},\chi=\mu|_{(A,e)\times(e,A)}$ if and only if the 2-cocycle
$\displaystyle \frac{\mu_0}{h{_{\mid A\times A}}}$ on $A\times A$ satisfies the conditions of Proposition \ref{lift} with respect to the action
$^{\overline t}(a,b)=(^t a,^{\phi(t^{-1})} b)$ of $S/A$ on $A\times A$, the 2-cocycle $\overline\rho=(\rho,\rho^{op})
\in Z^2(S/A,A\times A)$, and some $\psi\in Z^1(S/A,Fun(A\times A,\mathbb C^\times))$, $\nu\in C^2(T,\mathbb C^\times)$.

2) The condition $(\mu\circ(\vee\times\vee))^{-1}\cong \mu^{(e,x)}$, for some $x\in S^{op}$ such that $f^2=Ad(x)$, implies
the conditions $\mu^l(a,b)\cong (\mu^r(b^{-1},a^{-1}))^{-1},\chi(a,b)=\chi(b^{-1},a^{-1})$; in particular, the bicharacter $\chi$
is symmetric.
\end{corollary}

\begin{dm} 1) Lemma \ref{complot} and relation $(\omega\otimes\omega^{op})|_{L\times L\times L}=\partial^2 h$ allow to use
Remark \ref{extcohom}, 1) which gives the needed result. Note that
the relation between $\mu$ and extensions of $\mu^l,\mu^r$ and $\chi$ (in fact, of $\mu_0^l,\mu_0^r$ and $\chi_0) $ to $L$ is:
$$
\mu((a,b,t),(a',b',t'))=\mu^l((a,b,t),(a',b',t'))\mu^r((a',b',t),(a,b,t'))\times
$$
\begin{equation} \label{mu}
\times\chi((a,b,t),(a',b',t')),\ \text{for\ all}\ a,a',b,b'\in A,t,t'\in T.
\end{equation}

2) Direct computation.
\hfill\end{dm}

Let us summarize the above considerations. We will denote $\psi:=\mu_{T,A\times A}$, $\nu:=\mu_{T,T}$ (see Preliminaries),
$\varepsilon(M):=f^{-1}(M^{-1})$, for all $M\in S/A$, and $(BrPic(Vec_S^\omega))_{(A,\varepsilon)}$ the subset of $BrPic(Vec_S^\omega)$
composed by order two elements attached to $L<S\times S^{op}$, i.e., to the couple $(A,\varepsilon)$.

\begin{lemma} \label{vee}
Homomorphisms $c:\mathbb Z/2\mathbb Z\to\pi_1$ are indexed by collections $(A,\mu^r,\chi,\varepsilon,\psi,\nu)$, where $A\triangleleft S$ is abelian and equipped with a symmetric non-degenerate bicharacter $\chi$, $\mu^r\in (\Omega_{A,\omega})^S$, $\varepsilon\in Aut(S/A)$
such that $\varepsilon^2$ is inner, $\psi\in Z^1(S/A,Fun(A\times A,\mathbb C^\times))$, $\nu\in C^2(S/A,\mathbb C^\times)$ such that:

\begin{description}
\item [(i)] $\displaystyle\frac{\mu_0}{^{\overline t}\mu_0}=\partial^1_{A\times A}\psi(t,\cdot)$, for any $t\in S/A$,
where the $\omega\otimes\omega^{op}$-2-cocycle $\mu_0$ is defined above;

\item [(ii)]  $\partial^2\nu(t,t',t")=\psi(t,\overline\rho(t',t''))\displaystyle\frac{\mu_0(^t\overline\rho(t',t''),\overline\rho(t,t't''))}
{\mu_0(\overline\rho(t,t'),\overline\rho(tt',t''))},\ \forall t,t',t"\in S/A$;

\item [(iii)] There exist $k\in C^1(A\times A,\mathbb C^\times)$, $q\in C^1(S/A,\mathbb C^\times)$ such that:
\begin{align*}
\mu_0[ \mu_0\circ ( \vee \times \vee)]  &=  \partial^1 k,\\
\psi[\psi \circ( \epsilon \times \vee)] &=  \partial^0 \tilde k,\\
\nu [\nu \circ( \epsilon \times \epsilon)] &=  \displaystyle \frac{\partial^1 q}{k \circ \overline \rho},
\end{align*}
where $\tilde k$ is constant map $t \mapsto k$ on $S/A$;

\item [(iv)] $(BrPic(Vec_S^\omega))_{(A,\varepsilon)}$ is not empty.
\end{description}

Two collections, $(A,\mu^r,\chi,\varepsilon,\psi,\nu)$ and $(A',\mu'^r,\chi',\varepsilon',\psi',\nu')$, define the same homomorphism if and only if:
$A=A', \varepsilon=\varepsilon'$ in $Out(S/A)$, $\chi=\chi'$, $\mu^r\cong\mu'^r$ and there exist $\eta\in C^1(A\times A,\mathbb C^\times)$, $\phi \in C^1(S/A,\mathbb C^\times)$ such that $\displaystyle \frac{\mu_0'}{\mu_0} = \partial^1\eta
$, $\displaystyle \frac{\psi'}{\psi} = \partial^0\tilde \eta$, where $\tilde\eta$ is constant map $t \mapsto \eta$,
and $\displaystyle \frac{\nu'}{\nu} = \displaystyle \frac{\partial^1\phi}{\eta \circ \overline \rho}$.
\end{lemma}
\begin{dm}
1) Any homomorphism $c$ is defined by an equivalence class of invertible $c(0)$-bimodule categories $c(1)$ such that $c(1)\cong c(1)^{op}$, i.e.,
of the form $\mathcal M(L,\mu)$ (see Proposition \ref{M(L,mu) inv}), where $\mu\cong(\mu\circ(\vee\otimes\vee))^{-1}$. The subgroup $L$ of
$S\times S^{op}$ is of the form $L(A,f)$ - see Lemmas \ref{toutecrire} and \ref{a86}, this gives $A$, $\epsilon(t) := f^{-1}(t^{-1})$. Now we can use Corollary \ref{triple} which describes the relation between $\mu$ and triples $(\mu^l,\mu^r,\chi)$ in terms $\psi$ and $\nu$ and gives the relations $\mu^l(a,b)\cong (\mu^r(b^{-1},a^{-1}))^{-1},\chi(a,b)=\chi(b^{-1},a^{-1})$. To do this, due to condition (iv), we can use as $h$ a particular fixed $\overline\mu$ corresponding
to one of invertible $c(0)$-bimodule categories attached to the couple $(A,\varepsilon)$ and remark that $\overline\mu/\mu\in Z^2(L,\mathbb C^\times)$ can be chosen normal changing, if necessary, a representative $\mu$ in the same class. Finally, one can check, using Remark \ref{extcohom}, 2), that conditions (iii) are equivalent to the relation $\mu\cong(\mu\circ(\vee\otimes\vee))^{-1}$.

2) Two $c(0)$-bimodule categories, $\mathcal M(L(A,f),\mu)$ and $\mathcal M(L(A',f'),\mu')$, are equivalent if and only if $A= A'$ and there exists $x \in S$ such that $L'=Ad(1,x)L$ and $Ad(x)\circ f = f'$ for some $x\in S$
(this gives relations between $(A,\varepsilon)$ and $(A',\varepsilon')$ - see Lemma \ref{a86} and Remark \ref{cross2}) and $\mu'\cong \mu^{(e,x)}$. In order to show that this last condition is equivalent to the above conditions, one can apply Remark \ref{extcohom}, 2) to the coboundary $\displaystyle \frac{\mu'}{\mu^{(e,x)}}$ and its restrictions  $\displaystyle \frac{\mu'_0}{\mu_0}, \frac{\psi'}{\psi},\displaystyle \frac{\nu'}{\nu}$.
\hfill\end{dm}

\begin{remark} \label{omega=1}
When $\omega=1$, condition (iv) follows from other conditions because we can take in the proof $h=1$.
\end{remark}

\end{subsubsection}

\begin{subsubsection}
{\bf  Fusion ring structures for $\mathcal C=c(0)\oplus c(1)$} In order to equip the category $\mathcal C$ with a tensor product and a duality, we have to define an involutive
$c(0)$-bimodule equivalence $\gamma:c(1)\to c(1)^{op}\cong c(1)$ by $\gamma:M\mapsto M^*$, for any $M\in Irr(\mathcal M(L,\mu))$, such that
$(s\cdot M\cdot t)^*=t^{-1}\cdot M\cdot s^{-1}$. This implies some restrictions on the structure of
the fusion ring of $\mathcal C$ (see also \cite{Lip}, Lemma 2.3):

\begin{proposition}
\label{crabe01}
1) Given an invertible $c(0)$-bimodule category $c(1)=$ \newline $\mathcal M(L(A,f),\mu)$ equivalent to $c(1)^{op}$, let us suppose
that $c(0)\oplus c(1)$ is a fusion category with some tensor product $\otimes$. Then all possible $\mathbb Z/2\mathbb Z$-graded fusion rings
for $c(0)\oplus c(1)$ are characterized by elements  $\delta\in S/A$ such that $f^2 = Ad(\delta)$ and $f(\delta) = \delta ^{-1}$. Namely, for
all $s,s'\in S,M,N\in Irr(\mathcal M(L(A,f),\mu)))$:
$$
s^*= s^{-1},\ \ M^*= f(M)\delta,
$$
$$
s\otimes s'=  ss',\ \ s\otimes M = p(s) M, \ \ M\otimes s= Mf^{-1}(p(s)^{-1}),
M\otimes N^*=\underset{x \in MN^{-1}}\oplus x.
$$
2) Let $(f,\delta)$ and $(f',\delta')$ be as above, then  $c(0)$-bimodule categories $c(1)=\mathcal M(L(A,f),\mu))$ and $c(1)'=\mathcal
M(L(A,f'),\mu'))$ are equivalent if and only if there is $x\in S$, such that
$$
f'= Ad(xA)\circ f, \  \mu' = \mu \circ Ad(1,x)\   \rm{in} \  H^2(L(A,f'),\mathbb C^\times).
$$
3) $\mathbb Z/2\mathbb Z$-graded fusion rings for $c(0)\oplus c(1)$ and $c(0)\oplus c(1)'$ are isomorphic if and only if there is $F\in Aut(S)$ such that:
$$
F(A)=A,\ F(\delta) = \delta',\  F \circ f \circ F^{-1} = f'.
$$
\end{proposition}

\begin{dm} 1) Relations $s^*=s^{-1},s\otimes s'=ss'$ are obvious, relations
$s\otimes M = p(s) M, \ \ M\otimes s= Mf^{-1}(p(s)^{-1})$ follow from Lemmas \ref{rightaction} and \ref{toutecrire}.
The map $\gamma:M\mapsto M^*$ is involutive and, for any $N\in S/A$, $(M\cdot(au(N))^* = (au(N))^{-1}\cdot M^*$, so
$\gamma(Mf^{-1}(N^{-1})) =  N^{-1}\gamma(M)$. If $M=\bf 1$ and $N=f(N')^{-1}$, this gives $\gamma(N') = f(N')\gamma(\bf 1)$, and
if $N'=\delta:=\gamma(\bf 1)$, we have $f(\delta)= \delta^{-1}$.

The $\mathbb Z/2\mathbb Z$-grading implies that $M\otimes N=\underset{x\in S}\oplus n_x x$, where $n_x\in\mathbb N$, and the properties
of duality - that $n_e=1$ if and only if $M=N^*$ and $n_e=0$ if not. This gives $n_x=1$ if and
only if $M=xN^*$ (the set of such $x$ is an $A$-coset) and $n_x=0$ if not. In the above mentioned terms we have
\begin{align*}
M\otimes N
&= \underset{\{x \in S / M^*.x
=N\}}\oplus   x
&= \underset{\{ x \in S /  f(M)\delta f(p(x)^{-1}) = N\}}\oplus x\\
&=
\underset{x \in M\delta^{-1} f(N^{-1})} \oplus x
\end{align*}
which gives $M\otimes N^*=\underset{x \in MN^{-1}} \oplus x$,  so $N= {\bf 1} = A$ is the only simple object of $c(1)$ such that, for any other
simple object $M$ of $c(1)$, one has: $M \otimes N^* = \underset{x \in M} \oplus  x$, hence $\bf 1$ is an invariant of fusion rings.

2) Follows from Lemma \ref{a86}.

3) Let $(f,\delta)$ and $(f',\delta')$ be two pairs as above, then any isomorphism $\varphi$ of $\mathbb Z/2\mathbb Z$-graded fusion rings for $c(0)\oplus c(1)$
and $c(0)\oplus c(1)'$ is given by a pair $(F,\phi)$ where $F\in Aut(S)$ and $\phi : S/A \to S/A$ are such that, for any $s \in S,\ M,N \in S/A$,
one has:
$$
\phi(s\otimes M) = F(s) \otimes \phi(M), \phi(\bf 1)= \bf 1, \phi(M \otimes s)= \phi(M) \otimes F(s),
$$
$$
\phi(M^*)= \phi(M)^*,  \varphi(M \otimes N^*)= \phi(M) \otimes \phi(N)^*.
$$

The first equality gives $\phi(p(s)M) = p(F(s))\phi(M)$, so $F(A)\subset A$ and $F$ factors through $A$.
Let us write the last one can as $F(M)F(N)^{-1} = \phi(M)\phi(N)^{-1}$; together with $\phi(\bf 1)= \bf 1$,
this gives $\phi  = F$. The remaining equalities can be written as
$\phi(Mf^{-1}(p(s)^{-1}) = \phi(M)f'^{-1}(p(F(s^{-1}))$ and
$\phi(f(M)\delta) = f'(\phi(M))\delta'$, respectively, which can be resumed by:
\begin{equation}
\label{fusion2''}
f' = F \circ f \circ F^{-1},\ \delta' =  F(\delta).
\end{equation}

The converse also holds: given $(f,\delta)$ and $F\in Aut(S)$ such that $F(A) \subset A$, let us define $\phi:=F$ and $f'$, $\delta'$ by
(\ref{fusion2''}); then routine computations show that $f'$ is an anti-isomorphism of $S/A$ such that $f'^2 = Ad(\delta')$ and $f'(\delta')
= \delta'^{-1}$. If, moreover, there is $x\in S$ such that $f' = ad(xA)\circ f$ and if one defines  $\mu' := \mu \circ Ad(1,x)$,
then $\mathcal M(L(A,f'),\mu')$ is an invertible $c(0)$-bimodule category equivalent to its opposite and such that $\phi$ associated with
$(F,\phi)$ gives an isomorphism of the corresponding $\mathbb Z/2\mathbb Z$-graded fusion rings. \hfill \end{dm}
\end{subsubsection}
\end{subsection}

\begin{subsection} {$Vec_S^\omega$-bimodule equivalences $M_{g,h}$}

Given a group homomorphism $c:\mathbb Z/2\mathbb Z\to BrPic(c(0)): 0\mapsto\overline{c(0)},$ $1\mapsto\overline{c(1)}:=
\overline{\mathcal M(L,\mu)}$, there exist
(by definition of an invertible $c(0)$-bimodule category) $c(0)$-bimodule equivalences $M_{g,h}:c(g)\boxtimes_{c(0)}c(h)\to c(gh)\
(g,h\in\{0,1\})$. They are defined by linear functors respecting the
fusion rule of Proposition \ref{crabe01}: $M_{0,0}:t\boxtimes s\mapsto ts$, $M_{0,1}:t\boxtimes M\mapsto t\cdot M$, $M_{1,0}:M\boxtimes
t\mapsto M\cdot t$, $M_{1,1}:M\boxtimes N\mapsto \oplus_{M=xN^*}x$ equipped with natural
isomorphisms $M_{g,h}^l$ and $M_{g,h}^r$.

The $c(0)$-bimodule category structure on $c(0)$ is defined by $\mu^l(s,t,x)=\omega(s,t,x)$, $\mu^r(x,t,s)=\omega^{-1}(x,t,s), \chi(t,x,s)=
\omega(t,x,s)$, so the definition of a $c(0)$-bimodule functor and relations (\ref{lmod}), (\ref{rmod}), (\ref{leftbimod}), (\ref{rightbimod})
give, for all $t,x,y\in S,\ M\in S/A$:
$$
M_{0,0}^l(t,x\boxtimes y)=\omega(t,x,y)id_{txy},\ M_{0,0}^r(x\boxtimes y,t)=\omega^{-1}(x,y,t)id_{xyt},
$$
$$
M_{1,0}^l(t,M\boxtimes x)=\tilde\chi(t,M,x)id_{tMx},\ M_{1,0}^r(M\boxtimes x,t)=\tilde\mu^r(M,x,t)id_{Mxt},
$$
$$
M_{0,1}^l(t,x\boxtimes M)=\tilde\mu^l(t,x,M)id_{txM},\ M_{0,1}^r(x\boxtimes M,t)=\tilde\chi^{-1}(x,M,t)id_{xMt}.
$$
\begin{lemma} \label{balanced}
Given $c(1):=\mathcal M(L,\mu)$, there exists an equivalent $c(0)$-bimodule category for which
\begin{equation} \label{hyp1}
\tilde\mu^l(s,t,M)=\tilde\mu^r((stM)^*,s,t).
\end{equation}
In this equivalent category we have, up to a functor isomorphism:
\begin{equation} \label{l11}
M^l_{1,1}(t,M\boxtimes_{\mathcal D}N):=\oplus_{M=x\cdot N^*}[\tilde\mu^r(M^*t^{-1},t,x)]^{-1}id_{tx},
\end{equation}\begin{equation} \label{r11}
M^r_{1,1}(M\boxtimes_{\mathcal D}N,t):=\oplus_{M=x\cdot N^*}[\tilde\mu^r(M^*,x,t)]^{-1}id_{xt}.
\end{equation}
\end{lemma}

\begin{dm}
$M_{1,1}^l$ and $M_{1,1}^r$ must be defined by the equalities
\begin{align*}
M_{1,1}^l(t,M\boxtimes_{\mathcal D} N)=\underset {M=x\cdot N^*} \oplus\gamma^l(t,M,N,x)id_{tx},\\
M_{1,1}^r(M\boxtimes_{\mathcal D} N,t)= \underset {M=x\cdot N^*}\oplus \gamma^r(x,M,N,t)id_{xt},
\end{align*}
where the functions $\gamma^l(t,M,N,x)$ and $\gamma^r(x,M,N,t)$ satisfy, for all  $s,t \in S$,
$M,N\in Irr(\mathcal M)$, the conditions
$$
\omega(s,t,x)\gamma^l(st,M, N, x)
= \gamma^l(t,M,N,x)\gamma^l(s,tM,N,tx)\tilde \mu^l(s,t,M),
$$
$$
\omega^{-1}(x,s,t)\gamma^r(x,M,N,st)= \gamma^r(x,M,N,s)\gamma^r(xs,M,Ns,t)\tilde \mu^r(N,s,t),
$$
$$
\omega(s,x,t)\gamma^l(s,M,N,x) \gamma^r(sx,sM,N,t)=\gamma^r(x,M,N,t)\gamma^l(s,M,Nt,xt).
$$
The first of them with $M=N^*,x=e$ gives, denoting $g^l(t,M):=\gamma^l(t,M,M^*,e)$:
$$
\gamma^l(t,M,N,x)=\frac{g^l(tx,N^*)}{g^l(x,N^*)}[\tilde\mu^l(t,x,N^*)]^{-1},
$$
and, similarly, the second condition gives, denoting $g^r(t,M):=\gamma^r(e,M,M^*,t)$:
$$
\gamma^r(x,M,N,t)=\frac{g^r(xt,M)}{g^r(x,M)}[\tilde\mu^r(M^*,x,t)]^{-1}.
$$
Now the last condition becomes equivalent to the relation
$$
\tilde\mu^l(s,t,M)=\tilde\mu^r((stM)^*,s,t)\frac{\gamma(t,M)\gamma(s,t\cdot M)}{\gamma(st,M)},
$$
where $\gamma(t,M):= \displaystyle \frac{g^r(t,t\cdot M)}{g^l(t,M)}$, and one can check that $M_{1,1}^l$ and $M_{1,1}^r$
as above define $M_{1,1}$ as a surjective $c(0)$-bimodule functor.

Next, let us show that, passing to an isomorphic $c(0)$-bimodule functor, we can choose $g^l(t,M)\equiv 1$
and $g^r(t,M)=\gamma(t,t^{-1}\cdot M)$. Indeed, the definition of a $c(0)$-bimodule functor isomorphism
implies that two equivalences, $M_{1,1}$ and $M_{1,1}'$, corresponding to two couples, $(g^l, g^r)$ and
$(g'^l, g'^r)$, respectively, are isomorphic if and only if there is a natural transformation
$$
\alpha(M\boxtimes_{\mathcal D}N)=\oplus_{M=xN^*}\alpha(M,x,N)id_x,
$$
where $\alpha(M,x,N)$ satisfies the following system of two equations:
$$
\alpha(M,x,N)\frac{g^l(tx,N^*)}{g^l(x,N^*)}=\frac{g'^l(tx,N^*)}{g'^l(x,N^*)}\alpha(t\cdot M,tx,N),
$$
$$
\alpha(M,x,N)\frac{g^r(xt,M)}{g^r(x,M)}=\frac{g'^r(xt,M)}{g'^r(x,M)}\alpha(M,xt,N\cdot t).
$$
Putting $M=N^*$, $x=e$, and denoting $\delta(t,N^*):=\alpha(t\cdot N^*,t,N)$, we have:
$$
g^l(t,N^*)=\delta(t,N^*)g'^l(t,N^*),\ \ g^r(t,N^*)=\delta(t,(N\cdot t)^*)g'^r(t,N^*).
$$
As $\delta(t,N^*)$ can be arbitrary, we can choose $\delta(t,N^*)=g^l(t,N^*)$ and $g'^l(t,N^*)\equiv 1$ which gives
$g'^r(t,N^*)=\gamma(t,(N\cdot t)^*)$.

Finally, replacing $\tilde\mu^l(s,t,M)$ by $\tilde\mu^l(s,t,M)\displaystyle \frac{\gamma(t,M)\gamma(s,t\cdot M)}{\gamma(st,M)}$ and also
$\tilde\chi(s,M,t)$ by $\tilde\chi(s,M,t)\displaystyle \frac{\gamma(s,M\cdot t)}{\gamma(s,M)}$ and taking the same $\tilde\mu^r(M,s,t)$, we
pass to the equivalent $c(0)$-bimodule category, where the relations (\ref{hyp1}), (\ref{l11}) and (\ref{r11}) are true.\hfill\end{dm}

Let us deduce from (\ref{hyp1}) some properties of $\tilde\chi(s,M,t)\ (s,t\in S,M\in S/A)$.

\begin{remark} \label{restore} Theorem 3.4(iii) from \cite{Naidu} and relation (\ref{rightbimod}) imply the existence of
$\tilde\eta_{p(x)}\in C^1(S,C)$ (defined up to an element of $Z^1(S,C)$) such that map
$\tilde\chi(x,M,y)\tilde\eta_{p(x)}(M,y)\in Z^1(S,C),\tilde\eta_{\bf{1}}\equiv 1$. Similarly, (\ref{leftbimod}) and
(\ref{hyp1}) imply that, for any fixed $y\in S$,  one must have $\tilde\chi(x,M,y)\tilde\eta^{-1}_{p(y^{-1})}((xM)^*,x)\in\underline Z^1(S,C)$.
Then Shapiro's lemma shows that $\tilde\chi(x,M,y)$ is completely defined, up to a 1-coboundary, by one of its restrictions,
$\tilde\chi(x,{\bf 1},b)$ or $\tilde\chi(a,{\bf 1},y)$, where $a,b\in A$.
\end{remark}
\begin{lemma} \label{passation}
For any $x\in S$, there is $\eta_{p(x)}\in C^1(A,\mathbb C^\times)$ such that
$$
\tilde\chi(x,{\bf 1},b)=\chi(\kappa_{{\bf 1},x},b)\eta^{-1}_{p(x)}(b).
$$
\end{lemma}

\begin{dm}
Fixing $\tilde\eta_{p(x)}$ and applying the restriction map to (10), we have
$$
\frac{\eta'_{p(x)}(a)\eta'_{p(x)}(b)}{\eta'_{p(x)}(ab)}=\frac{\tilde\chi(x,{\bf 1},ab)}{\tilde\chi(x,{\bf 1},a)\tilde\chi(x,{\bf 1},b)},
$$
where $\eta'_{p(x)}(b):=\tilde\eta_{p(x)}(b,{\bf 1})$. So, for any $x\in S,b\in A$, $\tilde\chi(x,{\bf 1},b)\eta'_{p(x)}(b)=r_x(b)$ is
the unique character on $A$; this can be written as $\tilde\chi(u(p(x)),{\bf 1},b)\eta'_{p(x)}(b)$ $=\chi(\kappa_{{\bf 1},x},b)^{-1}r_x(b)$,
where the character in the right hand side depends only on $p(x)$. Since $\eta'_{p(x)}$ is defined modulo
$\hat A$, we have the needed equality, where $\eta_{p(x)}(b):=\eta'_{p(x)}(b)\chi(\kappa_{{\bf 1},x},b)r_x^{-1}(b)$.
\hfill \end{dm}
\begin{corollary} \label{simple}
(i) As $\tilde\chi(x,N,y)$ satisfies (9), then, putting there $N={\bf 1},M=p(y)$,  and $z=b$, we get
$$
\tilde\chi(x,M,b)=\frac{\chi(\kappa_{M,x},b)\eta_{M}(b)}{\eta_{x\cdot M}(b)}.
$$
Similar formula holds for $\tilde\chi(a,M,y)\ (a\in A,y\in S,M\in S/A)$. In particular, $\tilde\chi(a,M,b)=\chi(a,b)$, for all $a,b\in A$.

(ii) The following equalities hold, up to a product of two 1-coboundaries:
$$
\tilde \chi(a,M,x)=\tilde \chi(x,x^{-1}M^*,a),\
\tilde \chi(x,M,a)=\tilde \chi(a,M^*x^{-1},x)\ \text{for\ all}\ a\in A,x\in S.
$$
\end{corollary}

\begin{dm}
(i) follows from Lemma \ref{passation} and (ii) - from Lemma \ref{passation} and the symmetry of $\chi$.
\hfill\end{dm}

Rewriting (\ref{leftbimod}) and (\ref{rightbimod}) with the usage of (\ref{hyp1}), then comparing again with
these relations and, finally, taking in mind Corollary \ref{simple} (ii), we get:
\begin{corollary} \label{equiv}
If $\tilde\mu^l$ and $\tilde\mu^r$ satisfy (\ref{hyp1}), then the 2-cochain
\begin{equation} \label{hyp1'}
\beta(z,N,x):=\frac{\tilde\chi(x,x^{-1}N^*z^{-1},z)}{\tilde\chi(z,N,x)}
\end{equation}
is a product of  two 1-coboundaries.

Vice versa, let $\tilde\mu^r$
satisfy (\ref{rmod}) and $\tilde\chi$ satisfy the above mentioned properties. Then $\tilde\mu^l$ defined by (\ref{hyp1}) satisfies
(\ref{lmod}) and (\ref{leftbimod}).
\end{corollary}
\end{subsection}
\begin{subsection}{Vanishing of $O_3(c)$. Associativity isomorphisms. Main result}
\label{assoc}
\begin{subsubsection}
{\bf Quasi-tensor product  on $\mathcal C=c(0)\oplus c(1)$} Now, according to \cite{ENO2}, 8.4, we have to equip $\mathcal C$ with {\it a quasi-tensor product}, i.e., with a bifunctor
$\otimes:\mathcal C\times\mathcal C\to\mathcal C$ such that $\otimes\circ(\otimes\times Id_{\mathcal C})\cong\otimes\circ(Id_{\mathcal C}\times\otimes )$.
If we know a group homomorphism $c:\mathbb Z/2\mathbb Z\to\pi_1$ equipped with $c(0)$-bimodule equivalences $M_{g,h}:c(g)\boxtimes_{\mathcal D} c(h)\cong
c(gh)$, then \cite{ENO2}, Theorem 8.4 says that such a bifunctor exists if and only if the cohomological obstruction $O_3(c)$ vanishes or, in other
terms, if the following $c(0)$-bimodule autoequivalences are isomorphic to $Id$, for all $f,g,h\in \mathbb Z/2\mathbb Z$:
\begin{equation} \label{T}
T_{f,g,h}:M_{fg,h}(M_{f,g}\boxtimes_{c(0)} Id_{c(h)})(Id_{c(f)}\boxtimes_{c(0)}M_{g,h}^{-1})
M_{f,gh}^{-1}:c(fgh)\to c(fgh).
\end{equation}

Equivalently, we have to find conditions on parameters of $c$ such that the $c(0)$-bimodule functors $F_{f,g,h}:=M_{f,gh}\circ[id\boxtimes M_{g,h}]:
\C_f\boxtimes_{c(0)}\C_g\boxtimes_{c(0)} \C_h\to\C_{fgh}$ and $G_{f,g,h}:=M_{fg,h}\circ[M_{f,g}\boxtimes id]:
\C_f\boxtimes_{c(0)}\C_g\boxtimes_{c(0)} \C_h\to\C_{fgh}$ are isomorphic. In this case, they define the corresponding associativity isomorphisms:
\begin{equation} \label{alpha}
\alpha_{f,g,h}:M_{f,gh}(Id_{c(f)}\boxtimes_{c(0)}M_{g,h})\cong M_{fg,h}(M_{f,g}\boxtimes_{c(0)}Id_{c(h)}).
\end{equation}

First, we compute $F_{f,g,h}^l,F_{f,g,h}^r,G_{f,g,h}^l$ and $G_{f,g,h}^r$ using the
formulas for left and right bimodule structures of a composition of two bimodule functors:
$$
(F_2\circ F_1)^l(t,M)=F_2^l(t,F_1(M))\circ F_2(F_1^l(t,M)),
$$
\begin{equation}
(F_2\circ F_1)^r(M,t)=F_2^r(F_1(M),t)\circ F_2(F_1^r(M,t)).
\end{equation}

This and the fact that $[id\boxtimes M_{g,h}]^l$ and $[M_{f,g}\boxtimes id]^r$ are identities, give:
$$
F_{f,g,h}^l(t,X\boxtimes Y\boxtimes Z)=M^l_{f,gh}(t,X\boxtimes(Y\cdot Z)),
$$
$$
F_{f,g,h}^r(X\boxtimes Y\boxtimes Z,t)=M^r_{f,gh}(X\boxtimes(Y\cdot Z),t)M_{f,gh}[id_X\boxtimes M^r_{g,h}(Y\boxtimes Z,t)],
$$
$$
G_{f,g,h}^r(X\boxtimes Y\boxtimes Z,t)=M^r_{fg,h}((X\cdot Y)\boxtimes Z,t),
$$
$$
G_{f,g,h}^l(t,X\boxtimes Y\boxtimes Z)=M^l_{fg,h}(t,(X\cdot Y)\boxtimes Z)M_{fg,h}[M^l_{f,g}(t,X\boxtimes Y)\boxtimes id_Z].
$$
\end{subsubsection}
\begin{subsubsection}
{\bf Associativity isomorphisms for $\mathcal C=c(0)\oplus c(1)$} The associativity isomorphisms $\alpha_{f,g,h}(X,Y,Z)$, if they exist, 
must satisfy the system of relations:
$$
\alpha_{f,g,h}(t\cdot X,Y,Z)F^l_{f,g,h}(t,X\boxtimes Y\boxtimes Z)=\alpha_{f,g,h}(X,Y,Z)G_{f,g,h}^l(t,X\boxtimes Y\boxtimes Z),
$$
$$
\alpha_{f,g,h}(X,Y,Z\cdot t)F^r_{f,g,h}(X\boxtimes Y\boxtimes Z,t)=\alpha_{f,g,h}(X,Y,Z)G_{f,g,h}^r(X\boxtimes Y\boxtimes Z,t).
$$

At the same time, if this system has a solution, this means that the $c(0)$-bimodule functors $F_{f,g,h}$ and $G_{f,g,h}$ are isomorphic.
We compute step by step, using (\ref{hyp1}) and (\ref{hyp1'}):

1) $f=g=h=0$.
$F^l_{0,0,0}(t,x\boxtimes y\boxtimes z)=\omega(t,x,yz)id_{txyz},$
$F^r_{0,0,0}(x\boxtimes y\boxtimes z,t)=\omega(x,yz,t)\omega(y,z,t)id_{xyzt},$
$G^r_{0,0,0}(x\boxtimes y\boxtimes z,t)=\omega(xy,z,t)id_{xyzt},$
$G^l_{0,0,0}(t,x\boxtimes y\boxtimes z)=\omega(t,xy,z)\omega(t,x,y)id_{txyz}$,
$$
\alpha_{0,0,0}(x,y,z)=\omega(x,y,z)id_{xyz}.
$$
2) $f=1, g=h=0$.
$F^l_{1,0,0}(t,M\boxtimes x\boxtimes y)=\tilde\chi(t,M,xy)id_{tMxy},$
$F^r_{1,0,0}(M\boxtimes x\boxtimes y,t)=\tilde\mu^r(M,xy,t)\omega^{-1}(x,y,t)id_{Mxyt},$
$G^r_{1,0,0}(M\boxtimes x\boxtimes y,t)=\tilde\mu^r(Mx,y,t)id_{Mxyt}$,
$G^l_{1,0,0}(t,M\boxtimes x\boxtimes y)=\tilde\chi(t,M\cdot x,y)\tilde\chi(t,M,x)id_{tMxy}$,
$$
\alpha_{1,0,0}(M,x,y)=[\tilde\mu^r(M,x,y)]^{-1}id_{Mxy}.
$$
3) $f=g=0,h=1$.
$F^l_{0,0,1}(t,x\boxtimes y\boxtimes M)=\tilde\mu^l(t,x,y\cdot M)id_{txyM},$
$F^r_{0,0,1}(x\boxtimes y\boxtimes M,t)=
\tilde\chi^{-1}(x,y\cdot M,t)\tilde\chi^{-1}(y,M,t)id_{xyMt},$
$G^r_{0,0,1}(x\boxtimes y\boxtimes M,t)=\tilde\chi^{-1}(xy,M,t)\times$ $\times id_{xyMt}$,
$G^l_{0,0,1}(t,x\boxtimes y\boxtimes M)=\tilde\mu^l(t,xy,M)\omega(t,x,y)id_{txyM}$,
$$
\alpha_{0,0,1}(x,y,M)=\tilde\mu^l(x,y,M)id_{xyM}.
$$
4) $f=h=0,g=1$.
$F^l_{0,1,0}(t,x\boxtimes M\boxtimes y)=\tilde\mu^l(t,x,M\cdot y)id_{txMy},$
$F^r_{0,1,0}(x\boxtimes M\boxtimes y,t)=
\tilde\chi^{-1}(x,M\cdot y,t)\tilde\mu^r(M,y,t)id_{xMyt},$
$G^r_{0,1,0}(x\boxtimes M\boxtimes y,t)=\tilde\mu^r(x\cdot M,y,t)id_{xMyt}$,
$G^l_{0,1,0}(t,x\boxtimes M\boxtimes y)=
\tilde\mu^l(t,x,M)\tilde\chi(t,x\cdot M,y)id_{txMy}$,
$$
\alpha_{0,1,0}(x,M,y)=\tilde\chi(x,M,y)id_{xMy}.
$$
5) $f=g=1,h=0$.
$F^l_{1,1,0}(t,M\boxtimes N\boxtimes s)=\oplus_{M=xN^*}[\tilde\mu^l(t,xs,(N\cdot s)^*)]^{-1}id_{txs},$
$F^r_{1,1,0}(M\boxtimes N\boxtimes s,t)=\oplus_{M=xN^*}[\tilde\mu^r(M^*,xs,t)]^{-1}\tilde\mu^r(N,s,t)id_{xst}$, also one gets
$G^r_{1,1,0}(M\boxtimes N\boxtimes s,t)=\oplus_{M=xN^*}\omega^{-1}(x,s,t)id_{xst}$,
$G^l_{1,1,0}(t,M\boxtimes N\boxtimes s)=\oplus_{M=xN^*}\omega(t,x,s)
 \times[\tilde\mu^l(t,x,N^*)]^{-1}id_{txs}$
$$
\alpha_{1,1,0}(M,N,s)=\oplus_{M=xN^*}\tilde\mu^r(M^*,x,s)id_{xs}.
$$
6) $f=0, g=h=1$.
$F^l_{0,1,1}(t,s\boxtimes M\boxtimes N)=\oplus_{M=xN^*}\omega(t,s,x)id_{tsx},$
$F^r_{0,1,1}(s\boxtimes M\boxtimes N,t)=
\oplus_{M=xN^*}\omega^{-1}(s,x,t)[\tilde\mu^r(M^*,x,t)]^{-1}id_{sxt},$
$G^r_{0,1,1}(s\boxtimes M\boxtimes N,t)=\oplus_{M=xN^*}
[\tilde\mu^r((s\cdot M)^*,sx,t)]^{-1}id_{sxt},$
$G^l_{0,1,1}(t,s\boxtimes M\boxtimes N)=\oplus_{M=xN^*}[\tilde\mu^l(t,sx,N^*)]^{-1}$
$\times\tilde\mu^l(t,s,M)id_{tsx}$,
$$
\alpha_{0,1,1}(s,M,N)=\oplus_{M=xN^*}[\tilde\mu^l(s,x,N^*)]^{-1}id_{sx}.
$$
7) $f=h=1,g=0$.
$F^l_{1,0,1}(t,M\boxtimes s\boxtimes N)=\oplus_{Ms=xN^*}[\tilde\mu^l(t,x,(s\cdot N)^*)]^{-1}id_{tx},$
$F^r_{1,0,1}(M\boxtimes s\boxtimes N,t)=\tilde\chi^{-1}(s,N,t)\oplus_{Ms=xN^*}
[\tilde\mu^r(M^*,x,t)]^{-1}id_{xt}$, and also
$G^r_{1,0,1}(M\boxtimes s\boxtimes N,t)=\oplus_{Ms=xN^*}
[\tilde\mu^r((M\cdot s)^*,x,t)]^{-1}id_{xt}$, and finally  one gets
$G^l_{1,0,1}(t,M\boxtimes s\boxtimes N)=
\tilde\chi(t,M,s)\times\oplus_{Ms=xN^*}[\tilde\mu^l(t,x,N^*)]^{-1}id_{tx}$.

Then
$$
\alpha_{1,0,1}(M,s,N)=\oplus_{Ms=xN^*}\alpha(M,s,N,x)id_x,
$$

where  function $\alpha(M,s,N,x)$ satisfies the following two equations:
$$
\alpha(t\cdot M,s,N,tx)[\tilde\mu^l(t,x,(s\cdot N)^*)]^{-1}=
\alpha(M,s,N,x)[\tilde\mu^l(t,x,N^*)]^{-1}\tilde\chi(t,M,s).
$$
$$
\alpha(M,s,N\cdot t,xt)[\tilde\mu^r(M^*,x,t)]^{-1}\tilde\chi^{-1}(s,N,t)=
$$
$$
=\alpha(M,s,N,x)[\tilde\mu^r((M\cdot s)^*,x,t)]^{-1}.
$$

Putting $M=(s\cdot N)^*,x=e$ and denoting $\alpha(s,N):=\alpha((s\cdot N)^*,s,N,e)$, one can deduce,
for arbitrary $M,N$ and $x$ such that $M\cdot s=x\cdot N^*$:
$$
\alpha(M,s,N,x)=\alpha(s,N)\tilde\chi(x,x^{-1}\cdot M,s),
$$
$$
\alpha(M,s,N,x)=\alpha(s,N\cdot x^{-1})\tilde\chi(s,N\cdot x^{-1},x).
$$
Comparing these equalities and using (\ref{hyp1'}), where one can write
$$
\beta(z,N,s)=\frac{\alpha(z,N)}{\alpha(z,N\cdot s)}\ \ \text{with}\ \alpha(z,N)\in C^1(S,C),
$$
one has the following
\begin{lemma} \label{iso101}
The $c(0)$-bimodule autoequivalence $T_{1,0,1}$ is isomorphic to the identity, and we have:
\begin{equation} \label{101}
\alpha_{1,0,1}(M,s,N)=\alpha(s,N)\oplus_{Ms=xN^*}\tilde\chi(x,x^{-1}\cdot M,s)id_x.
\end{equation}
\end{lemma}

8) $f=g=h=1$.
$F^l_{1,1,1}(t,M\boxtimes N\boxtimes K)=\oplus_{N=yK^*}\tilde\chi(t,M,y)id_{tMy},$
$F^r_{1,1,1}(M\boxtimes N\boxtimes K,t)=\oplus_{N=yK^*}\frac{\tilde\mu^r(M,y,t)}
{\tilde\mu^r(N^*,y,t)}id_{Myt},$
$G^r_{1,1,1}(M\boxtimes N\boxtimes K,t)=\oplus_{M=xN^*}\tilde\chi(x,K,t)^{-1}id_{xKt},$
$G^l_{1,1,1}(t,M\boxtimes N\boxtimes K)=\oplus_{M=xN^*}\frac{\tilde\mu^l(t,x,K)}{\tilde\mu^l(t,x,N^*)}id_{txK}$.

As $\alpha_{1,1,1}(M,N,K):\oplus_{M=xN^*}x\cdot K\mapsto\oplus_{N=yK^*}M\cdot y$ (let us note that $x\cdot K=M\cdot y$),
this isomorphism is defined by an $|A|\times|A|$-matrix
$$
(\alpha(M,N,K,x,y)id_{x\cdot K})_{M=xN^*;N=yK^*},
$$
whose coefficients satisfy the system of the following relations
$$
\alpha(t\cdot M,N,K,tx,y)\tilde\chi(t,M,y)=\alpha(M,N,K,x,y)\frac{\tilde\mu^l(t,x,K)}{\tilde\mu^l(t,x,N^*)},
$$
$$
\alpha(M,N,K\cdot t,x,yt)\frac{\tilde\mu^r(M,y,t)}
{\tilde\mu^r(N^*,y,t)}=\alpha(M,N,K,x,y)\tilde\chi(x,K,t)^{-1}.
$$

Putting here $M=N^*=K,x=y=e$, we have that  $\alpha(t\cdot N^*,N,N^*,t,e)=\alpha(N^*,N,N^*,e,e)$ which we denote
by $\tau(N)$, $\alpha(N^*,N,N^*\cdot t,e,t)=\tau(N)$. Inserting these expressions again into
the above equations, we have
$$
\alpha(M,N,K,x,y)=\tau(N)\tilde\chi^{-1}(x,N^*,y).
$$

Thus, we have proved
\begin{lemma} \label{iso111}
The associativity isomorphism
$\alpha_{1,1,1}(M,N,K)$ is defined by  matrix
\begin{equation} \label{111}
(\tau(N)\tilde\chi^{-1}(x,N^*,y)id_{x\cdot K})_{M=xN^*;N=yK^*}.
\end{equation}
\end{lemma}
\begin{remark} \label{obstr}
The above results show that the obstruction $0_3(c)$ vanishes.
\end{remark}
\end{subsubsection}
\begin{subsubsection}
{\bf Vanishing of the obstruction $O_4(c,M)$} As $H^4(\mathbb Z/2\mathbb Z,\mathbb C^\times)=\{0\}$, the obstruction $O_4(c,M)$  (see \cite{ENO2}, 8.6) 
vanishes automatically, so there is a choice of $\alpha_{f,g,h}$ satisfying the pentagon equations (see \cite{ENO2}, (51)) which can be given a form 
similar to that in \cite{TY}:
$$
M_{f,gh,k}(id_f\boxtimes_{c(0)}\alpha_{g,h,k})\circ\alpha_{f,gh,k}(Id_{\mathcal C_f}\boxtimes_{c(0)}M_{g,h}\boxtimes_{\mathcal D}Id_{\mathcal C_k})
\circ M_{fgh,k}(\alpha_{f,g,h}\boxtimes_{c(0)} id_k)=
$$
\begin{equation} \label{pent}
=\alpha_{f,g,hk}(Id_{\mathcal C_f}\boxtimes_{c(0)}Id_{\mathcal C_g}\boxtimes_{c(0)}M_{h,k})\circ\alpha_{fg,h,k}(M_{f,g}\boxtimes_{c(0)}
Id_{\mathcal C_h}\boxtimes_{c(0)}Id_{\mathcal C_k}).
\end{equation}

Moreover, given two systems of associativity isomorphisms, $\alpha_{f,g,h}$ and $\alpha'_{f,g,h}$, \cite{ENO2}, Theorem 8.9 shows that $\beta(f,g,h)=\alpha_{f,g,h}(\alpha'_{f,g,h})^{-1}$ can be viewed as an element of $Z^3(\mathbb Z/2\mathbb Z,\mathbb C^\times)$, and that they give equivalent $\mathbb Z/2\mathbb Z$-extensions of $c(0)$ if and only if $\beta\in B^3(\mathbb Z/2\mathbb Z)$. As the only nontrivial 3-cocycle on $\mathbb Z/2\mathbb Z$ is defined by $\beta(1,1,1)=-1$, there are exactly 2 different classes of associativity isomorphisms which differ by a sign of $\tau(N)$. The simplest choice of representatives in these classes corresponds to $\alpha(t,M)\equiv 1$ and $\tau(N)\equiv\tau$. Exactly like in \cite{TY}, one can find easily from the pentagon equation for $\alpha_{1,1,1}$, using the non-degenaracy of $\chi$, that $\tau=\pm |A|^{-1/2}$.
\end{subsubsection}

\begin{subsubsection}
{\bf  Main result}
Let us recall, for the convenience of the reader, some notations, for other notations see Preliminaries. For a triple $(\mu^r,\mu^l,\chi)$, where $\mu^r\in(\Omega_{A,\omega})^S, \mu^l\cong(\mu^r\circ(\vee\times\vee))^{-1}, (x,y)^\vee:=(y^{-1},x^{-1}),\ \forall (x,y)\in (S\times S^{op})$, we denote $\mu_0((a,b),(c,d)):=\mu^l(a,c)\mu^r(b,d)\chi(a,d)$
- the $\omega\otimes\omega^{op}$-2-cocycle on $A\times A$; the action of $S/A$ on $A\times A$ is
denoted by $^{\overline t}(a,b):=(^{t}a,^{\varepsilon(t)} b)$, and we put $^{\overline t}\mu_0((a,b),(c,d)):=\mu_0(^{\overline t}(a,b),^{\overline t}(c,d))$. Finally, we introduce the 3-cochain on $S/A$:
$$
\zeta(t,t',t"):=\psi(t,\overline\rho(t',t''))\frac{\mu_0(^t\overline\rho(t',t''),\overline\rho(t,t't''))}
{\mu_0(\overline\rho(t,t'),\overline\rho(tt',t''))},
$$
where $\overline\rho=(\rho,\rho^{op})\in Z^2(S/A,A\times A)$ and $\rho(t,t'),\rho^{op}(t,t')=^{(t't)^{-1}}\rho (t',t)$ are 2-cocycles
giving, together with  the action $a\mapsto ^t a$ of $S/A$ on $A$, the structure of twisted semidirect product on $S$ and $S^{op}$,
respectively.

\begin{theorem} \label{main}  $\mathbb Z/2\mathbb Z$-extensions $\mathcal C$ of $c(0)=Vec^\omega_S$ are parameterized, up to a tensor equivalence, by collections $(A,\chi,\mu^r,\tau,\varepsilon,\delta,\psi,\nu)$, where:

- $A$ is an abelian normal subgroup of $S$ equipped with a symmetric non-degenerate bicharacter $\chi$, $\mu^r$ represents an equivalence class $(\Omega_{A,\omega})^S$ modulo restrictions on $A$ of $Z^2(S,\mathbb C^\times)$, $\tau=\pm|A|^{-1/2}$;

- $\varepsilon\in Aut(S/A)$, $\delta\in S/A$ such that $\varepsilon^2=Ad(\delta)$, $\varepsilon(\delta)=\delta$, and the set  $(BrPic(Vec_S^\omega))_{(A,\varepsilon)}$ is not empty;

- $\psi\in Z^1(S/A,Fun(A\times A,\mathbb C^\times))$ such that $\displaystyle \frac{\mu_0}{^{\overline t}\mu_0}=\partial^1 \psi,\ \text{for\ any}\ t\in S/A$, $\nu\in C^2(S/A,\mathbb C^\times)$ satisfying $\zeta=\partial^2\nu,
\mu_0[ \mu_0\circ ( \vee \times \vee)] = \partial^1 k,\
\psi[\psi \circ( \epsilon \times \vee)] =  \partial^0 \tilde k,\
\nu [\nu \circ( \epsilon \times \epsilon)] = \displaystyle \frac{\partial^1 q}{k \circ \overline \rho}$,
where $k\in C^1(A\times A,\mathbb C^\times),\ q\in C^1(S/A,\mathbb C^\times),\ \tilde k$ is constant map $t \mapsto k$ on $S/A$.

Namely, $\mathcal C=c(0)\oplus c(1)$ with $Irr(c(0))=S,Irr(c(1))=S/A$, the fusion rule is: $x^*=x^{-1},M^*=\varepsilon^{-1}(M^{-1})\delta,
x\otimes y=xy, x\otimes M=p(x)M,M\otimes x=M\varepsilon(p(x)), M\otimes N^*=\underset{x \in MN^{-1}}\oplus x\ (\text{where}\ x,y\in S,\ M,N\in S/A)$.

The associativity isomorphisms are defined, for any
$x,y,z \in S,\ K,L,M\in Irr(c(1))$, by 2-cochains $\tilde\mu^r$ and $\tilde\chi$ induced,respectively, from $\mu^r$
and from the extension of $\chi$ obtained by $\psi$ and $\nu$:
\begin{align*}
\alpha_{0,0,0}(x,y,z)
&= \omega(x,y,z)id_{xyz} \\
\alpha_{1,0,0}(K,x,y)
&= [\tilde\mu^r(K,x,y)]^{-1}id_{Kxy}\\
\alpha_{0,1,0}(x,K,y)
&= \tilde \chi(x,K,y)id_{xKy}\\
\alpha_{0,0,1}(x,y,K)
&= \tilde\mu^r((xyK)^*,x,y)id_{xyK}\\
\alpha_{0,1,1}(x,K,L)
&= \oplus_{K=sL^*} [\tilde\mu^r(K^*x^{-1},x,s)]^{-1}id_{xs}\\
\alpha_{1,1,0}(K,L,x))
&= \oplus_{K=sL^*} \tilde\mu^r(K^*,s,x)id_{sx}\\
\alpha_{1,0,1}(K,x,L)
&= \oplus_{Kx=sL^*} \tilde \chi(s,(xL)^*,x)id_{s} \\
\alpha_{1,1,1}(K,L,M)
\rm \ \ {is \ \ defined \ \ by\ the} \ \ matrix \ \ &(\tau\tilde \chi^{-1}(s,L^*,t)id_{sM})_{K=sL^*;L=tM^*}.
\end{align*}

Two collections, $(A,\chi,\mu^r,\tau,\varepsilon,\delta,\psi,\nu)$ and $(A',\chi',\mu'^r,\tau',\varepsilon',\delta',\psi',\nu')$, define equivalent
$\mathbb Z/2\mathbb Z$-extensions of $Vec^\omega_S$ if and only if:

-\hskip 0.5cm  $A=A',\tau=\tau'$ and there are $F\in Aut(S)$ and $\varphi\in C^2(S^{op},\mathbb C^\times)$
such that $\omega\circ (F\times F\times F)/\omega=\partial^2\varphi,\ F(A)=A$
(so $F$ factors through $A$), $F(\delta)=\delta'$ and, modulo inner automorphisms: $\mu'^r\cong\varphi|_{A\times A}\cdot\mu^r\circ(F\times F),\chi'=\chi\circ(F\times F),\ F\circ \varepsilon=\varepsilon'\circ F$;

-\hskip 0.5cm  $\displaystyle \frac{\mu_0\circ(F\times F)}{\mu'_0}=\partial^1\eta$, $\displaystyle \frac{\psi\circ F} {\psi'}=\partial^0\tilde\eta$, $\displaystyle \frac{\nu\circ F}{\nu'}=\displaystyle \frac{\partial^1(\phi\circ F)}
{\eta\circ\overline\rho}$, where $\eta\in C^1(A\times A,\mathbb C^\times)$, $\tilde\eta:t\mapsto\eta$ is constant function
on $S/A$, and $\phi\in C^1(S/A,\mathbb C^\times)$.

\end{theorem}

\begin{dm} (i) If $\mathcal C=c(0)\oplus c(1)$ is a fusion category, then $c(1)\cong c(1)^{op}$ is an invertible $c(0)$-bimodule category,
so it corresponds to a collection $(A,\mu^r,\chi,\varepsilon,\psi,\nu)$ described in Lemma \ref{vee}. The existence of $\delta$ and
its properties follow from Proposition \ref{crabe01}, the existence and the value of $\tau$, as explained above, - from \cite{ENO2},
Theorem 8.9 and from the pentagon equation for $\alpha_{1,1,1}$. Thus, we have associated with $\mathcal C$ the collection $(A,\mu^r,\chi,\tau,\varepsilon,\delta,\psi,\nu)$.

(ii)
Vice versa, given a collection $(A,\chi,\mu^r,\tau,\varepsilon,\delta,\psi,\nu)$ as above, Lemma \ref{vee} allows to construct an
invertible $c(0)$-bimodule category $c(1)$ such that $c(1)\cong c(1)^{op}$, and Proposition \ref{crabe01} - a fusion ring. Now, Lemmas \ref{iso101},
\ref{iso111} and the triviality of $H^4(\mathbb Z/2\mathbb Z)$ show that there is a tensor product on $\mathcal C=c(0)\oplus c(1)$
giving a structure of a fusion category, with two choices of $\tau$.

(iii) Let $(A,\chi,\mu^r,$ $\tau,\varepsilon,\delta,\psi,\nu)$ and $(A',\chi',\mu'^r,\tau',\varepsilon',\delta',\psi',\nu')$ be two
collections corresponding to equivalent $\mathbb Z/2\mathbb Z$-extensions, $\mathcal C$ and $\mathcal C'$, respectively. By definition, a tensor
equivalence of $\mathcal C$ and $\mathcal C'$ contains:

1) a tensor autoequivalence of $c(0)$ defined by a couple $(F,\varphi)$, where $F\in Aut(S)$ and $\varphi\in C^2(S,\mathbb C^\times)$
are such that $\omega\circ(F\times F\times F)/\omega=\partial^2\varphi$ (it is more convenient for us to pass to $c(0)^{op}$ and so to
deal with $\varphi\in C^2(S^{op},\mathbb C^\times)$);

2) an equivalence of $c(0)$-bimodule categories $c(1)=\mathcal M(L(A,f),\mu)$ and $c(1)'=\mathcal M(L(A',f'),\mu')$ which
 implies the following equalities with some $f^l,f^r\in C^1(S,C)$ (see Remark \ref{cross2}):
$$
\frac{\tilde\chi'^s}{\tilde\chi\circ(F\times F)}=\frac{\partial^0 f^l}{{\underline{\partial}^0 f^r}},\ \ \
\frac{(\tilde\mu'^r)^s}{\tilde\mu^r\circ(F\times F)}=\varphi\cdot\partial^1 f^r
$$
whose restriction to $A$ gives $(\chi')^s=\chi\circ(F\times F)$ and $\varphi|_{A\times A}\cdot\mu^r\circ(F\times F)\cong(\mu'^r)^s$
(these equalities become simpler when considered modulo inner automorphisms:
$\mu'^r\cong\varphi|_{A\times A}\cdot\mu^r\circ(F\times F),\chi'=\chi\circ(F\times F)$);

Lemma \ref{vee} gives the remaining list of relations between the components of the above collections.

3) an isomorphism of their fusion rings, so Proposition \ref{crabe01} implies that $F(A)=A$, $F(\delta)=\delta'$, and
$F\circ\varepsilon=\varepsilon'\circ F$ in $Out(S/A)$.

Finally, the explicit formulas for the associativity isomorphisms for two given $\mathbb Z/2\mathbb Z$-extensions of $c(0)$ and the above relations
between the two corresponding collections allow to construct an equivalence of these categories.
\hfill\end{dm}

\begin{corollary} \label{semi-direct}
If $S\cong A\rtimes(S/A)$ is a usual semidirect product, i.e., $\rho=1$, one can omit $\mu^r$ in the above parameterization
and the conditions on $\psi$ and $\nu$ are simpler: $\nu\in Z^2(S/A,\mathbb C^\times),\ \psi'\circ(\varepsilon\otimes\vee)=
(\psi\circ F)^{-1},\ \nu'\circ(\varepsilon\otimes\varepsilon)\cong(\nu\circ F)^{-1}$.
\end{corollary}
\begin{dm}
Let us fix a $\mathbb Z/2\mathbb Z$-extension $c(0)\oplus c(1)$ of $c(0)$ and the $\omega$-2-cocycle $\mu_1^r$ giving a structure of right $c(0)$-bimodule category on
$c(1)$. Let $\mathcal C$ be any other $\mathbb Z/2\mathbb Z$-extension with its $\mu_2^r$. We want to show, using tensor equivalences of the form $(Id,\varphi),
\ \varphi\in Z^2(S^{op},\mathbb C^\times)$, that there is a $\mathbb Z/2\mathbb Z$-extension
$\mathcal C'$ equivalent to $\mathcal C$ for which $\mu'^r_2=\mu^r_1$. We have, on the one hand, $\mu'^r_2=\varphi|_{A\times A}\mu^r_2$, and on the
other hand, $\mu^r_2=Z\cdot \mu^r_1$, where $Z\in Z^2(A,\mathbb C^\times)$. So, it suffices to choose as $\varphi$ any extension of $Z^{-1}$, if it exists.

But the conditions imposed on $(\mu^r_1,\chi_1)$ and $(\mu^r_2,\chi_2)$ imply that $Z=\displaystyle \frac{\mu^r_2}{\mu^r_1}$ satisfies all
the conditions of \cite{Kar}, Lemma 2.2.4, so its extension exists.

Thus, if $S\cong A\rtimes(S/A)$, one can omit $\mu^r$ in the above parameterization of $\mathbb Z/2\mathbb Z$-extensions of $c(0)$.
Also, as $\mu'^r_2=\mu^r_1$, we have $\eta=1$ in the conditions of Theorem \ref{main}, so they take the above mentioned simpler form.
\hfill\end{dm}

\begin{remark}
\label{rank 1}
a) If $\omega=1$, the condition $(BrPic(Vec_S^\omega))_{(A,\varepsilon)}$ is not empty follows from other conditions (see Remark \ref{omega=1});
if also $S\cong A\rtimes(S/A)$, one can choose $\mu^r=1$.

b) If $A=S$, Theorem \ref{main} gives the result of Tambara and Yamagami \cite{TY}.
This result was also obtained by the methods of graded fusion categories in \cite{ENO2}, Proposition 9.3, Example 9.4.

c) The case $A=\{e\}$ was treated in much more generality in \cite{ENO2}.
\end{remark}
\end{subsubsection}
\end{subsection}
\end{section}
\begin{section}{Examples}

\begin{subsection}{Non-isomorphic fusion rings if $S/A\cong \mathbb Z/p\mathbb Z,\ p$ is prime}

\begin{lemma} \label{cyclicrules}
If $S/A\cong \mathbb Z/p\mathbb Z,\ p$ is prime, and let  $n$ be the number of isomorphism classes of $\mathbb Z/2\mathbb Z$-graded
fusion rings with basis $X = S \cup S/A$, then:

(i) if $p= 2$,  one has: $n=2$;

(ii) if $\rho=1$, i.e $S\cong A\underset{\alpha}\rtimes \mathbb Z/p\mathbb Z$ is a usual semidirect product, where $\alpha$ is an action of $\mathbb Z/p\mathbb Z$ on $A$, then $n=3$ if $S$ is
abelian and $n=p+1$ otherwise.
\end{lemma}

\begin{dm}
(i) If $p=2$, then Proposition \ref{crabe01} gives $f=Id$ and $\delta=0$ or $\delta=1$ give two non-isomorphic fusion rings corresponding to the single subgroup
$L=L^{op}$ of $S\times S^{op}$.

(ii) If $p>2$, $S/A\cong \mathbb Z/p\mathbb Z = \{0,1,..,p-1\}$, then $f_1=id$ and $f_{-1}:x \mapsto -x$ are the only involutive elements of $Aut(S/A)$ which give two subgroups, $L_1=L_1^{op}$ and $L_{-1}=L_{-1}^{op}$ - see Proposition \ref{crabe01}.

Now, only $\delta = 0$ satisfies equation $f_1(\delta)=-\delta$; on the contrary, any $\delta\in S/A$ satisfies equation $f_{-1}(\delta)= -\delta$.
As $F\circ id\circ F^{-1}=id$, for any $F\in Aut(S/A)$, the fusion ring corresponding to the couple $(id,0)$ is not isomorphic to any other. Similarly, as $F(0)=0$, for any $F\in Aut(S/A)$, the fusion ring corresponding to the couple $(f_1,0)$ is not isomorphic to any other. Let us also remark that
$F\circ f_{-1}\circ F^{-1}=f_{-1}$, for any $F\in Aut(S/A)$.

Next, Proposition \ref{crabe01} shows that the only condition under which two couples, $(f_{-1},\delta)$ and $(f_{-1},\delta')\ (\delta,\delta'\neq 0)$, generate isomorphic fusion ring, is the existence of $F\in Aut(S/A)$ such that $F(\delta) = \delta'$ and $F$ can be extended to $G \in Aut(S)$ such that $G(A) = A$. If $S= A\underset\alpha\rtimes\mathbb Z/p\mathbb Z $, this is equivalent to
$$
\alpha_{\delta'\delta^{-1}x}G(a) = G(\alpha_{x}(a)),\ \text{for\ all}\ a \in A,x \in \mathbb Z/n\mathbb Z.
$$

Putting here $x=1$, we get $\alpha_\delta= \alpha_{\delta'}$, i.e, $\delta'\delta^{-1}\in Ker(\alpha)$, where $\alpha: \delta\to\alpha_\delta$
is a homomorphism from $S/A$ to $Aut(A)$. Conversely, if $\delta'\delta^{-1}\in Ker(\alpha)$, then one can put $G(a) = a$, for any $a \in A$, in order to extend $F$. As $p$ is prime, $Ker(\alpha)$ is $S$ when  $S$ is abelian or $\{e\}$ otherwise, and the result follows.
\hfill\end{dm}
\end{subsection}

\begin{subsection}{Number of non-equivalent $\mathbb Z/2\mathbb Z$-extensions in examples}

\begin{definition}
If $S$ is a finite group, $\omega\in H^3(S,\mathbb C^\times)$, $A\triangleleft S$ is abelian and
such that $\omega|_{A\times A\times A}=1$ in $H^3(A,\mathbb C^\times)$, we denote by $n(S,\omega,A)$ the number of
$\mathbb Z/2\mathbb Z$-graded fusion categories (up to equivalence) associated with it.
\end{definition}

\begin{proposition} \label{egypte}
If $S$ is abelian, $|S|=2p\ (p\ \text{is\ prime})$ and $A$ is a nontrivial subgroup of $S$, then $n(S,1,A)$ equals to 6 when $S\cong\mathbb Z/2p\mathbb Z,
A\cong \mathbb Z/2\mathbb Z,\ p>2$, to 8 when $S\cong\mathbb Z/4\mathbb Z$ or $S\cong \mathbb Z/2p\mathbb Z, A\cong \mathbb Z/p\mathbb
Z,\ p>2$, and to 16 when $S\cong(\mathbb Z/2\mathbb Z)^2$, for any of 3 its nontrivial subgroups.
\end{proposition}

\begin{dm}
If $S\neq\mathbb Z/4\mathbb Z$, all the subgroups $L$ are direct products of $A\times A$ and $S/A$.
If $|A|=2$, there is only one symmetric non-degenerate bicharacter
on $A$; if $A\cong\mathbb Z/p\mathbb Z,\ p>2$, \cite{W} tells us that there are 2 orbits in the set of non-degenerate symmetric bicharacters
on $A$ with respect to $Aut(A)$. As $S$ is abelian, $Inn(A)=\{Id\}$, as $S/A$ is cyclic, $H^2(S/A,\mathbb C^\times)=0$.

We have $H^1(S/A,A\times A)=0$ when $S\cong \mathbb Z/2p\mathbb Z,\ p>2$ because $|S/A|$ and $|A|$ are relatively prime
- see \cite{Kar}, Proposition 1.3.1.
According to Lemma \ref{cyclicrules}, the number of non isomorphic fusion rings equals to 2 if $|S/A|=2$ and to 3 otherwise.
As in all cases there are 2 choices for $\tau$, we already have the needed result for $S\cong \mathbb Z/2p\mathbb Z,\ p>2$.

To complete the proof if $S\cong(\mathbb Z/2\mathbb Z)^2$, it suffices to remark that $H^1(\mathbb Z/2\mathbb Z,$
$(\mathbb Z/2\mathbb Z)^2)$, with trivial action, is $(\mathbb Z/2\mathbb Z)^2$.

If $S\cong\mathbb Z/4\mathbb Z,\ A\cong\mathbb Z/2\mathbb Z$, then $S\cong A\underset \rho\rtimes S/A$ with trivial action of $S/A$ on $A$
and the 2-cocycle $\rho:S/A\to A$ given by $\rho(1,1)=1$, and $\rho=0$ otherwise. The subgroup $L<S\times S^{op}$ is isomorphic to
$(A\times A)\underset{\tilde \rho}\rtimes S/A$ with trivial action and the 2-cocycle $\tilde \rho:S/A\to
A\times A$ coming from $\rho$. As above, $H^1(S/A,A\times A)=(\mathbb Z/2\mathbb Z)^2$; in order to check condition 3) of Proposition \ref{lift},
remark that the function $\gamma(s,t,u):=\alpha_{(0,S/A),(A\times A,0)}(s,\tilde \rho(t,u))$ equals identically to 1 for two of them and to the
nontrival 3-cocycle on $\mathbb Z/2\mathbb Z$ for two others.

Now, the equation $\partial^2\alpha=\gamma$ has no solutions when $\gamma$ is a nontrival 3-cocycle because the left-hand side is clearly
a 3-coboundary. For two elements of $H^1(S/A,A\times A)$ for which the solution of the above equation exists, it is unique because $H^2(\mathbb Z/2\mathbb Z,\mathbb C^\times)=0$. Thus, taking in mind 2 choices for $\tau$, we see that the number of non-equivalent $\mathbb Z/2\mathbb Z$-extensions
in this example is 4 for any of the two non-isomorphic fusion rings.\hfill\end{dm}

\begin{remark} \label{bonderson}
It was shown earlier in \cite{Bo} (see also \cite{Lip}) that the number of fusion categories coming from $S\cong \mathbb Z/4\mathbb Z,\ \omega=1,\
A\cong \mathbb Z/2\mathbb Z,\ \delta=1$ is 4.
\end{remark}

\begin{proposition}
\label{odeur}
$n(D_p,1,\mathbb Z/p\mathbb Z)=8$, where $D_p$ is dihedral group, $p>2$ is prime, $A=\mathbb Z/p\mathbb Z$.
\end{proposition}

\begin{dm}
$D_p\cong A\rtimes\mathbb Z/2\mathbb Z$ where $A = \mathbb Z/p\mathbb Z$, the subgroup $L=L^{op}<D_p\times D_p^{op}$ is unique and isomorphic to
$(A\times A)\rtimes\mathbb Z/2\mathbb Z$ with the action $(a,b)\mapsto(-a,-b)$, so any bicharacter on $A\times A$ is $S/A$-invariant. Since $[S:A]=2$,
there are 2 non isomorphic fusion rings - see Lemma \ref{cyclicrules}.

As above, there are two orbits in the set of non-degenerate symmetric bicharacters on $A$ with respect to automorphisms of $S$ satisfying
the conditions of Theorem \ref{main}.

Finally, $H^2(S/A,\mathbb C^\times)=H^2(\mathbb Z/2\mathbb Z,\mathbb C^\times)=1$ and $H^1(S/A,A\times A)=1$ because $|S/A|$ and $|A\times A|$ are relatively prime
-see \cite{Kar}, Proposition 1.3.1.

Thus, non-equivalent $\mathbb Z/2\mathbb Z$-graded extensions of $Vec^1_{D_p}$ correspond to 2 choices of fusion rings, to 2 choices of $\chi$, and to 2 choices of $\tau$.
\hfill\end{dm}

\begin{proposition}
\label{A4}
$n(A_4,1,\mathbb Z/2\mathbb Z\times \mathbb Z/2\mathbb Z) = 8$, where $A=\mathbb Z/2\mathbb Z\times\mathbb Z/2Z$.
\end{proposition}

\begin{dm}
$A_4\cong A\rtimes\mathbb Z/3\mathbb Z$, where $A = \mathbb Z/2\mathbb Z \times \mathbb Z/2\mathbb Z$, and both subgroups, $L_1$ and $L_{-1}$, are isomorphic to
usual semidirect products of  type  $A^2\rtimes\mathbb Z/3\mathbb Z$. Since $[S:A]=3$, there are 4 non isomorphic fusion rings - see Lemma \ref{cyclicrules},
so we have to show that, for any fixed $\varepsilon,\delta$, there are exactly 2 fusion categories.

There are 4 symmetric non-degenerate bicharacters on $A:\ \chi_K((a_0,a_1),$ $(b_0,b_1)) = (-1)^{\Sigma^1_{i,j=0} k_{ij}a_ib_j}$, where $a_i,b_j,k_{ij}\in \{0,1\}$, $k_{ij}$ are matrix coefficients of a symmetric invertible $2\times 2$-matrix $K$.

Due to Theorem \ref{main}, in order to find non-equivalent $\mathbb Z/2\mathbb Z$-graded extensions of $Vec^1_{A_4}$, we have to look at the automorphisms of $A_4$ under which $A$ is stable. It is straightforward to show that inner automorphisms of $A_4$ generated by $Z/3Z$ transform one into another
the three bicharacters $\chi_K$ with $K\neq \begin{pmatrix} 0 & 1\\1 & 0\end{pmatrix}$ and that the fourth one is stable under $Aut(A)$. A simple computation gives that in the case of $L_1$ the three bicharacters in the same orbit are $S/A$-cohomologically invariant but not the fourth one; conversely, in the case of $L_{-1}$ only the fourth one is $S/A$-cohomologically invariant.

Finally, $H^2(S/A)=H^2(\mathbb Z/3\mathbb Z,\mathbb C^\times)=1$ and $H^1(S/A,A\times A)=1$ because $|S/A|$ and $|A\times A|$ are relatively prime - see \cite{Kar}, Proposition 1.3.1.

Thus, non-equivalent $\mathbb Z/2\mathbb Z$-graded extensions of $Vec^1_{A_4}$ correspond to 4 choices of fusion rings, to 1 choice of $\chi$, and to 2 choices of $\tau$.
\hfill\end{dm}
\end{subsection}
\end{section}


\bibliographystyle{ams-alpha}

\end{document}